\documentclass[10pt]{article}
\usepackage{amsmath,amsfonts,amsthm,amsopn}

\renewcommand{\ge}{\varepsilon} \newcommand{\R}{\mathbb{R}}
\newcommand{\C}{\mathbb{C}}

\DeclareMathOperator{\re}{Re}
\DeclareMathOperator{\im}{Im}

\newtheorem{theorem}{Theorem}[section]
\newtheorem{lemma}[theorem]{Lemma}

\newtheorem{proposition}[theorem]{Proposition}

\theoremstyle{definition}
\newtheorem{remark}[theorem]{Remark}

\begin{document}

\title{\textbf{Single--peaks for a magnetic Schr\"{o}dinger equation with critical growth}}

\author{Sara Barile\footnote{Email: \texttt{s.barile@dm.uniba.it}. Supported by MIUR, national project
    \textit{Variational and topological methods in the study of
      nonlinear phenomena}.} \\ Dipartimento di Matematica,
  Universit\`a di Bari \\ via Orabona 4, I-70125 Bari
  \and Silvia Cingolani\footnote{Email: \texttt{cingolan@poliba.it}. Supported by
    MIUR, national project \textit{Variational and topological methods
      in the study of nonlinear phenomena}.} \\ Politecnico di Bari \\ via Amendola 126/B, I-70126 Bari
  \and Simone
  Secchi\footnote{Email: \texttt{secchi@mat.unimi.it}. Supported by MIUR, national project
    \textit{Variational methods and nonlinear differential
      equations}.} \\ Dipartimento di Matematica, Universit\`a di
  Milano \\ via Saldini 50, I-20133 Milano.}
\maketitle

\begin{abstract}
  We prove existence results of complex-valued solutions for a
  semilinear Schr\"o\-dinger equation with critical growth under the
  perturbation of an external electromagnetic field.  Solutions are
  found via an abstract perturbation result in critical point theory,
  developed  in \cite{ab,ABC,am05}.
\end{abstract}

\bigskip


\bigskip \bigskip \bigskip

\section{Introduction}

\bigskip

This paper deals with some classes of elliptic equations which are
perturbation of the time-dependent nonlinear Schr\"{o}dinger
equation
\begin{equation}
  \frac{\partial \psi}{\partial t} = -\hbar^2 \Delta  \psi -
  |\psi|^{p-1} \psi
\end{equation}
under the effect of a magnetic field $B_\varepsilon$ and an electric
field $E_\varepsilon$ whose sources are small in $L^\infty$ sense.
Precisely we will study the existence of wave functions  $\psi
\colon \R^N \times \R\to \C$ satisfying  the nonlinear Schr\"odinger
equation
\begin{equation}\label{eq:ev}
  \frac{\partial \psi}{\partial t} = \left( \frac{\hbar}{i} \nabla
    -A_\varepsilon(x)\right)^2 \psi + W_\varepsilon(x)  \psi -
  |\psi|^{p-1} \psi
\end{equation}
where $A_\varepsilon (x)$ and $W_\varepsilon(x)$ are respectively a
magnetic potential and an electric one, depending on a positive
small parameter $\varepsilon >0$.  In the work, we assume that
$A_\varepsilon (x) = \varepsilon \ A(x)$, $W_\varepsilon(x)= V_0 +
\varepsilon^\alpha V(x)$, being $A \colon \R^N \to \R^N$ and  $V_0
\in \R$, $V \colon \R^N \to \R$, $\alpha \in [1,2]$.

On the right hand side of $(\ref{eq:ev})$ the operator $\left(
\frac{\hbar}{i} \nabla - A_\varepsilon \right)^2$ denotes the formal
scalar product of the operator $\frac{\hbar}{i} \nabla -
A_\varepsilon$ by itself, {\it i.e.}
\[
\left( \frac{\hbar}{i} \nabla -A_\varepsilon(x) \right)^2 \psi :=
-\hbar^2 \Delta \psi -\frac{2\hbar}{i}  A_\varepsilon  \cdot \nabla
\psi +  |A_\varepsilon|^2 \psi - \frac{\hbar}{i} \psi
\operatorname{div}A_\varepsilon
\]
being $i^2=-1$, $\hbar$ the Planck constant.

This model arises in several branches of physics, e.g. in the
description of the Bose--Einstein condensates and in nonlinear
optics (see \cite{avhesi,best,helffermorame,lupan}).

If $A$ is seen as the 1--form
\[
A = \sum_{j=1}^N A_j dx_i,
\]
then
\[
B_\varepsilon = \varepsilon \, dA = \varepsilon \sum_{j<k} B_{jk} \,
dx_j \wedge dx_k, \quad \text{where $B_{jk} =
  \partial_j A_k - \partial_k A_j$},
\]
represents the external magnetic field having source in $\varepsilon
A$ (cf. \cite{sulem}), while $E_\varepsilon = \varepsilon^\alpha
\nabla V(x)$ is the electric field. The fixed $\hbar >0$ the
spectral theory of the operator has been studied in detail,
particularly by Avron, Herbt, Simon \cite{avhesi} and Helffer
\cite{hel1,RS}.

The search of standing waves of the type
$\psi_\varepsilon(t,x)=e^{-i V_0 \hbar^{-1}t} u_\varepsilon(x)$
leads to find a complex-valued solution  $u \colon \R^N \to \C$ of
the semilinear Schr\"{o}dinger equation
\begin{equation}\label{tp}
\left( \frac{\hbar}{i} \nabla - \varepsilon A(x)\right)^2 u
+\varepsilon^\alpha V(x) u = |u|^{p-1} u \quad\hbox{in $\R^N$}.
\end{equation}
From a mathematical viewpoint, this equation has been studied in
several papers in the \emph{subcritical case} $1<p< (N+2)/(N-2)$. In
the pioneering paper \cite{EL}, M. Esteban and P.L. Lions proved the
existence of standing wave solutions to $(\ref{eq:ev})$ in the case
$V=1$ identically, $\varepsilon >0$ fixed, by a constrained
minimization. Recently variational techniques are been employed to
study equation $(\ref{tp})$ in the semiclassical limit $(\hbar \to
0^+)$.  We refer to \cite{ci,cs,ku,ss}. Recent results on
multi-bumps solutions are obtained in \cite{caotang} for bounded
vector potentials and in \cite{cs05} without any
$L^\infty$--restriction on $|A|$.

In the \emph{critical case}  $p= (N+2)/(N-2)$, we mention the paper
\cite{ariolisz} by Arioli and Szulkin where the potentials $A$ and
$V$ are assumed to be periodic, $\varepsilon >0$ fixed. The
existence of a solution is proved whenever $0 \notin \sigma \left(
\left( \frac{\nabla}{i}-A \right)^2 +V \right)$. We also cite the
recent paper \cite{chsz} by Chabrowski and Skulzin, dealing with
entire solutions of $(\ref{tp})$.

In the present paper we are concerned with the critical case $p=
(N+2)/(N-2)$, but $V$ and $A$ are not in general periodic
potentials.

When the problem is nonmagnetic and static, i.e. $A=0$, $V =0$, and
$\hbar =1$ then problem $(\ref{tp})$ reduces to the equation
\begin{equation} \label{eq:9}
-\Delta u = u^{\frac{N+2}{N-2}}, \quad \hbox{$u \in
D^{1,2}(\R^N,\C)$}.
\end{equation}
In Section 2 we prove that the least energy solutions to
$(\ref{eq:9})$ are given by the functions  $z= e^{i \sigma} z_{\mu,
\xi}(x)$, where
 \begin{equation} \label{rich}
z_{\mu, \xi}(x)=
 \kappa_N
\frac{\mu^{(\frac{N}{2}-1)}}{(\mu^2+|x-\xi|^2)^\frac{N-2}{2}}, \quad
\kappa_N=\left(N(N-2)\right)^\frac{N-2}{4}
\end{equation}
and they correspond to the extremals of the Sobolev imbedding
$D^{1,2}(\R^N, \C) \subset L^{2^*} (\R^N, \C)$ (cf. Lemma
\ref{exxx}).

The perturbation of $(\ref{eq:9})$ due to the action of an external
magnetic potential $A$ leads us to seek for complex--valued
solutions. In general, the lack of compactness due to the critical
growth of the nonlinear term produces several difficulties in facing
the problem by global variational methods. We will attack
$(\ref{tp})$ by means of a perturbation method in Critical Point
Theory, see \cite{ab,am99,am05}, and we prove the existence of a
solution $u_\ge$ to $(\ref{tp})$ that is close for $\ge$ small
enough to a solution of \eqref{eq:9}. After an appropriate finite
dimensional reduction, we find that stable critical points on $]0, +
\infty[ \times \R^{N+1}$ of a suitable functional $\Gamma$
correspond to points on ${Z}=\left\{e^{i \sigma}z_{\mu,\xi} \, : \,
  \sigma \in S^1, \, \mu > 0, \, \xi \in \R^N \, \right\} $ from which
there bifurcate solutions to $(\ref{tp})$ for $\varepsilon \ne 0$.
If $V$ changes its sign, we find at least two solutions to
$(\ref{tp})$. The main result of the paper is Theorem
\ref{mainresult}, stated in Section 5.

We quote the papers \cite{aap,c,cipi}, dealing with perturbed
semilinear equations with critical growth without magnetic potential
$A$.

\begin{remark} \label{rem:orbits} It is apparent that the compact
  group $S\sp 1$ acts on the space of solutions to $(\ref{tp})$. For
  simplicity, we will talk about solutions, rather than orbits of
  solutions.
\end{remark}

\noindent \textbf{Notation.} The complex conjugate of any number
$z\in\C$ will be denoted by $\bar z$. The real part of a number
$z\in\C$ will be denoted by $\re z$. The ordinary inner product
between two vectors $a,b\in {\R^N}$ will be denoted by $a \cdot b$.
We use the Landau symbols. For example $O(\ge)$ is a generic
function such that $\limsup_{\ge\to 0} \frac{O(\ge)}{\ge} < \infty$,
and $o(\ge)$ is a function such that $\lim\limits_{\ge\to 0}
\frac{o(\ge)}{\ge}=0$.  We will denote $D^{1,2}(\R^N,\C) =
\left\lbrace u \in L^{2^*}(\R^N,\C) \mid \int_{\R^N}
  |\nabla u|^2 \, dx < \infty \right\rbrace $, with a similar definition for
$D^{1,2}(\R^N,\R)$.


\section{The limiting problem}

\bigskip

Before proceeding, we recall some known facts about a couple of auxiliary
problems. Recall that $2^* = 2N/(N-2)$.

\medskip

\noindent (\textbullet) \ The problem
\begin{equation} \label{scalar}
\begin{cases}
  -\Delta u = |u|^{2^*-2}u &\hbox{in $\R^N$} \\
  u \in D^{1,2}(\R^N,\R).
\end{cases}
\end{equation}
possesses a smooth manifold of least-energy solutions
\begin{equation}
\widetilde{Z}=\left\{z_{\mu,\xi}=\mu^{-\frac{(N-2)}{2}} z_0
(\tfrac{x-\xi}{\mu}) \mid  \mu > 0, \, \xi \in \R^N   \right\}
\end{equation}
where
\begin{equation}
z_0(x)=\kappa_N \frac{1}{(1+|x|^2)^\frac{N-2}{2}}, \quad
\kappa_N=\left(N(N-2)\right)^\frac{N-2}{4}.
\end{equation}
Explicitly,
\begin{equation} \label{eq:10}
z_{\mu, \xi}(x)= \kappa_N \mu^{-\frac{(N-2)}{2}}
\frac{1}{\left(1+\left|\frac{x-\xi}{\mu} \right|^2 \right)^\frac{N-2}{2}}
= \kappa_N  \frac{\mu^{(\frac{N}{2}-1)}}{(\mu^2+|x-\xi|^2)^\frac{N-2}{2}}.
\end{equation}
These solutions are critical points of the Euler functional
\begin{equation} \label{euler}
\widetilde{f}_{0}(u)=\frac{1}{2}\int_{\R^N} |\nabla u |^{2}-\frac{1}{2^{*}}\int_{\R^N} |u_{+}|^{2^{*}} dx,
\end{equation}
defined on $D^{1,2}(\R^N,\R)\subset E$, and the following
\textit{nondegeneracy} property holds:
\begin{equation}
\ker \widetilde{f}_0'' (z_{\mu,\xi}) = T_{z_{\mu,\xi}} \widetilde{Z} \quad \hbox{for all $\mu>0$, $\xi \in \R^N$}.
\end{equation}

\medskip

\noindent (\textbullet\textbullet) \ Similarly, $f_0 \in
C^2(D^{1,2}(\R^N,\C))$ possesses a finite--dimensional manifold $Z$ of
least-energy critical points, given by
\begin{equation}
Z=\left\{e^{i \sigma} z_{\mu , \xi}:\sigma \in S^1, \mu > 0, \xi \in \R^N
\right\} \cong S^1 \times (0, + \infty) \times \R^N.
\end{equation}
More precisely, following the ideas of \cite{ku} and \cite{ss}, we
give the following characterization.
\begin{lemma}\label{exxx}
Any least-energy solution to the problem
\begin{equation}\label{complex}
\begin{cases}
  -\Delta u = |u|^{2^*-2}u &\text{in $\R^N$} \\
  u \in D^{1,2}(\R^N,\C)
\end{cases}
\end{equation}
is of the form $u=e^{i\sigma} z_{\mu,\xi}$ for some suitable $\sigma \in
[0,2\pi]$, $\mu>0$ and $\xi\in\R^N$.
\end{lemma}
\begin{proof}
  It is convenient to divide the proof into two steps.\\
  \textbf{Step 1:} Let $z_0=U$ the least energy solution associated to the
  energy functional (\ref{euler}) on the manifold
\begin{equation*}
M_{0,r}=\left\{v \in\, D^{1,2}(\R^N, \R) \setminus \{0\} \mid \int_{\R^N} {\left| \nabla v\right|}^2 \, dx=
\int_{\R^N} {\left|  v\right|}^{2^* }\, dx \right\}.
\end{equation*}
It is well-known that $z_0=U$ is radially symmetric and unique (up
to translation and dilation) positive solution to the equation
\eqref{scalar}. Let $b_{0,r}=b_r=\widetilde{f_0}(U)=
\widetilde{f_0}(z_0)$. In a similar way, we define the class
\begin{equation*}
M_{0,c}=\left\{v \in\, E \setminus \{0\} \mid \int_{\R^N} {\left| \nabla v\right|}^2 \, dx=
\int_{\R^N} {\left|  v\right|}^{2^* }\, dx \right\}
\end{equation*}
and denote by $b_{0,c}=b_c=f_0(v)$ on $M_{0,c}$. Let $\sigma \,\in\,\R$,
$\xi\,\in\,\R^N$, $\mu > 0$ ,$ \widetilde{v}(x)= z_{\mu, \xi}(x)$ positive
solution to (\ref{scalar}) and $\widetilde{U}= e^{i \sigma}\widetilde{v}=e^{i
  \sigma}z_{\mu, \xi}$ (i.e. $z_{\mu, \xi}= |\widetilde{U}(x)|$). It results
that
$\widetilde{U}=e^{i \sigma}z_{\mu, \xi}$ is a non-trivial least energy solution for $b_{0,c}=f_0(v)$ with $v\,\in \,M_{0,c}$.\\
\\
\textbf{Step 2:} The following facts hold:
 \begin{itemize}
 \item[(i)] $b_{0,c}= b_{0,r};$
 \item[(ii)] If $U_c= \widetilde{U}$ is a least energy solution of problem
   (\ref{complex}), then
\begin{equation*}
\left| \nabla |U_c |(x) \right|=\left| \nabla U_c (x) \right| \, \quad\text{and}\, \,   \re \left(i \overline{U_c}(x) \nabla U_c(x) \right)=0 \quad\text{for a.e. $x\,\in\,\R^N$.}
\end{equation*}
\item[(iii)] There exist $\sigma\,\in\,\R$ and a least energy solution $u_r :
  \R^N \to \R$ of problem (\ref{scalar}) with
\begin{equation*}
U_c(x)= e^{i \sigma} u_r(x)   \quad\text{ for a.e. $x\,\in\,\R^N$}
\end{equation*}
or, equivalently, the least energy solution $U_c$ for $b_{0,c}$ is the
following
 \begin{equation*}
U_c(x)= e^{i \sigma} u_r(x)= e^{i \sigma} z_{\mu,\xi}(x)    \quad\text{ for a.e. $x\,\in\,\R^N$}.
\end{equation*}
\end{itemize}
 Observe that
 \begin{equation*}
 b_{0,r}= \min_{v\,\in\,M_{0,r}} \widetilde{f_0}(v) \quad\text{and \, \, }   b_{0,c}= \min_{v\,\in\,M_{0,c}} f_0(v)
 \end{equation*}
 where $M_{0,r}$ and $M_{0,c}$ are the real and complex Nehari manifolds for
 $\widetilde{f_0}$ and $f_0$,
 \begin{eqnarray*}
M_{0,r}&=& \left\{v \in\, D^{1,2}(\R^N, \R) \setminus \{0\} \mid  {\widetilde{f}_0}'(v)[v]=0\right\}\\
 &=& \left\{v \in\, D^{1,2}(\R^N, \R) \setminus \{0\} \mid \int_{\R^N} {\left| \nabla v\right|}^2 \, dx=
\int_{\R^N} {\left|  v\right|}^{2^* }\, dx \right\}
 \end{eqnarray*}
 and
  \begin{eqnarray*}
M_{0,c}&=& \left\{v  \in  E \setminus \{0\} \mid  f'_0(v)[v]=0\right\}\\
 &=& \left\{v  \in  E \setminus \{0\} \mid \int_{\R^N} {\left| \nabla v\right|}^2 \, dx=
\int_{\R^N} {\left|  v\right|}^{2^* }\, dx \right\}
 \end{eqnarray*}
 So (i) is equivalent to
 \begin{eqnarray*}
 b_{0,r} &=& \min_{v\,\in\,M_{0,r}} \widetilde{f_0}(v)=\widetilde{f_0}(u_r)\\
 b_{0,c} &=& \min_{v\,\in\,M_{0,c}} f_0(v)=f_0(U_c)
 \end{eqnarray*}

\noindent \textit{Proof of} (i)--(iii).
Let $u\,\in\,E$ be given. For the sake of convenience, we introduce the
functionals
\begin{eqnarray*}
T(u)&=& \int_{\R^N} {\left| \nabla u\right|}^2 \, dx\\
P(u)&=& \frac{1}{2^*}\int_{\R^N} {\left|  u\right|}^{2^* }\, dx
\end{eqnarray*}
(resp. $\widetilde{T}(u)$ and $\widetilde{P}(u)$ as $u\,\in\, D^{1,2}(\R^N,
\R)$) such that
$f_0(u)=\frac{1}{2} T(u)- P(u)$ as $u\,\in\,E$ (resp. $\widetilde{f_0}(u)=\frac{1}{2} \widetilde{T}(u)- \widetilde{P}(u)$ as $u\,\in\, D^{1,2}(\R^N, \R)$).\\
Consider the following minimization problems
\begin{eqnarray*}
\sigma_r &=& \min \left\{\widetilde{T}(u) \mid u\,\in\, D^{1,2}(\R^N, \R), \widetilde{P}(u)=1  \right\}\\
\sigma_c &=& \min \left\{T(u) \mid u\,\in\, E, P(u)=1  \right\}
\end{eqnarray*}
Note that, obviously, there holds $\sigma_c \,\leq\,\sigma_r$. If we denote by
$u_*$ the Schwarz symmetric rearrangement (see \cite{blio}) of the positive
real valued function $|u| \,\in\,D^{1,2}(\R^N, \R)$, then, Cavalieri's
principle yields
\[
\int_{\R^N} {\left| u_*\right|}^{2^* }\, dx=\int_{\R^N} {\left| u\right|}^{2^*
}\, dx
\]
which entails $\widetilde{P}(u_*)=P(|u|)$. Moreover, by the
Polya-Sz\"{e}g\"{o} inequality, we have
\[
\widetilde{T}(u_*)=\int_{\R^N} {\left| \nabla u_*\right|}^{2}\, dx \leq
\int_{\R^N} {\left| \nabla |u| \right|}^{2}\, dx \leq \int_{\R^N} {\left|
    \nabla u\right|}^{2}\, dx=T(u)
\]
where the second inequality follows from the following diamagnetic inequality
\[
\int_{\R^N} {\left| \nabla |u| \right|}^{2}\, dx \leq \int_{\R^N} {\left|
    D^{\varepsilon} |u| \right|}^{2}\, dx \quad\text{for all $u\,\in\,
  H^{\varepsilon}_{A,V}$}
\]
with $ D^{\varepsilon}=\frac{\nabla}{i}-\varepsilon A$ and $A=0$. Therefore,
one can compute $\sigma_c$ by minimizing over the subclass of positive,
radially symmetric and radially decreasing functions $u \,\in\,D^{1,2}(\R^N,
\R)$. As a consequence, we have $ \sigma_r \, \leq \,\sigma_c$. In conclusion,
$ \sigma_r = \sigma_c$. Observe now that
\begin{eqnarray*}
b_{0,r}&=& \min \left\{\widetilde{f_0}(u) \mid \text{$u \in D^{1,2}(\R^N, \R) \setminus \left\{0\right\}$     is a solution to (\ref{scalar})}                      \right\},\\
b_{0,c}&=& \min \left\{f_0(u) \mid \text{$u \in E \setminus \left\{0\right\}$ is a solution to (\ref{complex})}                      \right\}.
\end{eqnarray*}
The above inequalities hold since any nontrivial real (resp. complex) solution
of (\ref{scalar}) (resp. (\ref{complex})) belongs to $M_{0,r}$ (resp.
$M_{0,c}$) and, conversely, any solution of $ b_{0,r}$ (resp. $ b_{0,c}$)
produces a nontrivial solution of (\ref{scalar}) (resp. (\ref{complex})).
Moreover, it follows from an easy adaptation of \cite[Th. 3]{blio} that $
b_{0,r}=\sigma_r $ as well as $ b_{0,c}=\sigma_c $. In conclusion, there holds
\[
 b_{0,r}=\sigma_r = b_{0,c}=\sigma_c
\]
which proves (i).

To prove (ii), let $U_c \colon \R^N \to\C$ be a least energy solution to
problem (\ref{complex}) and assume by contradiction that
\[
\mathcal{L}^N \left(\,\left\{x\,\in\,\R^N: \left| \nabla |U_c| \right|\,<
    \,\left| \nabla U_c \right| \right\}\, \right)\,>\, 0
\]
where $\mathcal{L}^N$ is the Lebesgue measure in $\R^N$. Then, we would get
$\widetilde{P}(|U_c|)=P(U_c)$ and
\[
\widetilde{P}(|U_c|)=\frac{1}{2^*} \int_{\R^N} {|U_c|}^{2^*}\,
dx= \frac{1}{2^*} \int_{\R^N} {\left| U_c \right|}^{2^*}\, dx =P(U_c)
\]
and
\[
\sigma_r \,\leq\, \int_{\R^N} {\left| \nabla |U_c| \right|}^{2}\, dx <
\int_{\R^N} {\left| \nabla U_c \right|}^{2}\, dx =\sigma_c
\]
which is a contradiction. The second assertion in (ii) follows by direct
computations. Indeed, a.e. in $\R^N$, we have
\[
\left| \nabla |U_c| \right|=\left| \nabla U_c \right| \quad\text{if and only
  if } \re U_c \left(\nabla \im U_c \right)= \im U_c \nabla \left( \re U_c
\right).
\]
If this last condition holds, in turn, a.e. in $\R^N$, we have
\[
\overline{U_c} \nabla U_c =\re U_c \nabla \left( \re U_c \right)+ \im U_c
\nabla \left( \im U_c \right)
\]
which implies the desired assertion.

Finally, the representation formula of (iii) $U_c(x)= e^{i \sigma} u_r(x)$ is
an immediate consequence of (ii), since one obtains $U_c= e^{i \sigma} |U_c|$
for some $\sigma \,\in\,\R$.
\end{proof}

\begin{remark}
For the reader's convenience, we write here the second derivative of $f_0$ at
any $z \in Z$:
\begin{multline}
\left\langle f''_{0}(z) v, w\right\rangle_{E}= \re \int_{\R^N} \nabla v
\cdot \, \overline{\nabla w} \, dx- \re \int_{\R^N} |z|^{2^{*}-2} v \overline
w \, dx \\ - \re (2^{*}-2) \int_{\R^N} |z|^{2^{*}-4} \re (z \overline v)z
\overline w \, dx.
\end{multline}
In particular, $f''_0(z)$ can be identified with a compact perturbation of
the identity operator.
\end{remark}

We now come to the most delicate requirement of the perturbation method.
\begin{lemma}
For each $z=e^{i \sigma} z_{\mu, \xi} \in Z$, there holds
\[
T_{z}Z= \ker f''_{0}(z) \qquad \hbox{for all $z \in Z$},
\]
where
\begin{equation}\label{tang}
T_{e^{i \sigma} z_{\mu, \xi}}Z=
\operatorname{span}_{\R}\left\{\frac{\partial e^{i \sigma} z_{\mu, \xi} }{\partial \xi_{1}},\dots, \frac{\partial e^{i \sigma} z_{\mu, \xi} }{\partial \xi_{N}},    \frac{\partial e^{i \sigma} z_{\mu, \xi} }{\partial \mu}, \frac{\partial e^{i \sigma} z_{\mu, \xi} }{\partial \sigma}= i e^{i \sigma} z_{\mu, \xi}   \right\}.
\end{equation}
\end{lemma}
\begin{proof}
The inclusion $T_{z}Z \subset \ker f''_{0}(z)$ is always true, see
\cite{ab}.  Conversely, we
prove that for any $\varphi \in \ker f''_{0}(z)$ there exist numbers $a_1, \dots, a_N$, $b$, $d \in \R$ such
that
\begin{equation} \label{eq:3}
\varphi=\sum_{j=1}^{N} a_j \frac{\partial e^{i \sigma} z_{\mu, \xi}
}{\partial \xi_{j}}+ b \frac{\partial e^{i \sigma} z_{\mu, \xi}
}{\partial \mu} + d i e^{i \sigma} z_{\mu, \xi}.
\end{equation}
If we can prove the following representation formul\ae, then
\eqref{eq:3} will follow.
\begin{eqnarray}
\re (\varphi  \overline{e^{ i \sigma}}) &=& \sum_{j=1}^{N} a_j  \frac{\partial  z_{\mu, \xi} }{\partial \xi_{j}}+ b \frac{\partial  z_{\mu, \xi} }{\partial \mu} \label{eq:re} \\
\im(\varphi  \overline{e^{ i \sigma}}) &=& d   z_{\mu, \xi}. \label{eq:im}
\end{eqnarray}
We will use a well-known result for the scalar case:
\begin{equation*}
\ker \widetilde{f}''_{0}( z_{\mu, \xi}) \equiv T_{z_{\mu, \xi}}\widetilde{Z}=
\operatorname{span}_{\R}\left\{\frac{\partial  z_{\mu, \xi} }{\partial \xi_{1}},\dots, \frac{\partial  z_{\mu, \xi} }{\partial \xi_{N}},    \frac{\partial  z_{\mu, \xi} }{\partial \mu}  \right\}
\end{equation*}
\textbf{Step 1: proof of \eqref{eq:re}}.
We wish to prove that $\re (\varphi \overline{e^{ i \sigma}}) \in \ker
\widetilde{f}''_{0}( z_{\mu, \xi})$.
Recall that $\varphi \in \ker f''_0(e^{ i \sigma} z_{\mu, \xi})$, so
\begin{equation} \label{eq:7}
\langle f''_0(e^{ i \sigma} z_{\mu, \xi})\varphi, \psi \rangle=0
\quad\hbox{for all $\psi \in E$}.
\end{equation}
Select $\psi=e^{ i \sigma} v$, with $v \in C_0^\infty(
\R^N,\R)$.
\begin{multline*}
0=\langle f''_0(e^{ i \sigma} z_{\mu, \xi})\varphi, v e^{i \sigma} \rangle
=\re \int
\nabla (\varphi e^{ -i \sigma}) \nabla v \\ -(2^* -2) \int_{\R^N}
|z_{\mu,\xi}|^{2^*-2} \re (\overline{e^{i \sigma}} \varphi)v - \int_{\R^N}
|z_{\mu,\xi}|^{2^*-2} \re (\overline{e^{i \sigma}} \varphi)v\\
= \int_{\R^N} \nabla (\re (\varphi \overline{ e^{ i \sigma}}) \nabla v -(2^*
-1) \int_{\R^N} |z_{\mu,\xi}|^{2^*-2} \re (\overline{e^{i \sigma}}
\varphi)v 
= \langle \widetilde{f}''_0 ( z_{\mu, \xi}) \re (\varphi
\overline{e^{i \sigma}}),v \rangle.
\end{multline*}
This implies that
\begin{equation*}
\re (\overline{e^{i \sigma}} \varphi) \in \ker \widetilde{f''_0}(
z_{\mu, \xi})\equiv T_{z_{\mu, \xi}} \widetilde{Z}
\end{equation*}
from which it follows
\begin{equation*}
\re (\varphi \overline{e^{ i \sigma}})=\sum_{j=1}^{N} a_j  \frac{\partial  z_{\mu, \xi} }{\partial \xi_{j}}+ b \frac{\partial  z_{\mu, \xi} }{\partial \mu}
\end{equation*}
for some real constants $a_1 ,\dots, a_N$ and $b$.

\textbf{Step 2: proof of \eqref{eq:im}.}  Test (\ref{eq:7}) on $\psi =
i e^{i \sigma}w \in E$ with $w:\R^N \to \R$.  We get
\begin{eqnarray}
  0 &=& \langle f''_0(e^{ i \sigma} z_{\mu, \xi}) \varphi , i e^{i \sigma}
  w \rangle =\re \int_{\R^N}  \nabla (-i \varphi e^{-i  \sigma}) \cdot \nabla w - \re
  \int_{\R^N}
 |z_{\mu, \xi}|^{2^* -2} (-i \varphi e^{-i \sigma}) w \notag \\
  &&
  \quad \hbox{[being $\re (- i \varphi e^{-i \sigma})= \im (\varphi
  e^{-i \sigma}) $]} \notag \\
  &=& \int_{\R^N} \nabla (\im (\varphi e^{-i \sigma}))
  \cdot \nabla w - \int_{\R^N} |z_{\mu, \xi}|^{2^*-2} \im (\varphi e^{-i
  \sigma}) w \notag \\
  &=& \int_{\R^N} \nabla (\im (\varphi \overline{ e^{i
  \sigma}})) \cdot \nabla w - \int_{\R^N} |z_{\mu, \xi}|^{2^*-2} \left[\im
  (\varphi \overline{e^{i \sigma}})\right]_{+} w.  \label{eq:11}
\end{eqnarray}
We can take $\mu = 1$ and $\xi = 0$, otherwise we perform the change of
variable $ x \mapsto \mu x + \xi$.

From \eqref{eq:11} we get that $u:=\im (\varphi \overline{e^{i
\sigma}})$ satisfies the equation
\begin{equation} \label{eq:8}
- \Delta u = \frac{N(N-2)}{(1+|x|^2)^2}
u \quad \hbox{in $D^{-1,2} (\R^N,\R)$}.
\end{equation}
We will study this linear equation by an inverse stereographic projections
onto the sphere $S^N$. Precisely, for each point $\xi \in S^N$, denote by $x$
its corresponding point under the stereographic projection $\pi$ from $S^N$ to
$\R^N$, sending the north pole on $S^N$ to $\infty$. That is, suppose $\xi =
(\xi_1,\xi_2,\dots,\xi_{N+1})$ is a point in $S^N$, $x=(x_1,\ldots,x_N)$, then
$\xi_i = \frac{2x_i}{1+|x|^2}$ for $1 \leq i \leq N$; $\xi_{N+1} =
\frac{|x|^2-1}{|x|^2+1}$.

Recall that, on a Riemannian manifold $(M,g)$, the conformal
Laplacian is defined by
\[
L_g = -\Delta_g + \frac{N-2}{4(N-1)} S_g,
\]
where $-\Delta_{g}$ is the
Laplace--Beltrami operator on $M$ and $S_g$ is the scalar curvature of $(M,g)$. It is known that
\[
L_g (\Phi (u)) = \varphi^{-\frac{N+2}{N-2}} L_{\delta} (u),
\]
where $\delta$ is the euclidean metric of $\R^N$, $\varphi (x) = \left(
  \frac{2}{1+|x|^2} \right)^{(N-2)/2}$ and

\[
\Phi \colon {D}^{1,2}(\R^N) \to H^1(S^n),\qquad
                     \Phi(u)(x)={u(\pi(x))\over \varphi(\pi(x))}
\]
is an isomorphism between $H^1(S^n)$ and $E:=D^{1,2}(\R^N)$.
Therefore, if $U = \Phi (u)$, then
\eqref{eq:8} changes into the equation
\begin{equation} \label{eq:23}
-\Delta_{g_0} U + \frac{N-2}{4(N-1)} S_{g_0} U = \frac{N(N-2)}{4} U,
\end{equation}
where $g_0$ is the standard riemannian metric on $S^N$, and
\[
S_{g_0} = N(N-1)
\]
is the constant scalar curvature of $(S^N,g_0)$. As a consequence, \eqref{eq:23} implies that
\[
-\Delta_{g_0}U =0,
\]
i.e. $U$ is an eigenfunction of $-\Delta_{g_0}$ corresponding to the
eigenvalue $\lambda =0$. But the point spectrum of $-\Delta_{g_0}$ is
completely known (see \cite{bs,bgm}), consisting of the numbers
\[
\lambda_k = k(k+N-1), \qquad k= 0,1,2,\dots
\]
with associated eigenspaces of dimension
\[
\frac{(N+k-2)! \, (N+2k-1)}{k! \, (N-1)!}.
\]
Hence we deduce that $k=0$, and $U$ belongs to an eigenspace of dimension
$1$. Since $z_{\mu,\xi}$ is a solution to \eqref{eq:8}, we conclude that
there exists $d\in\R$ such that
\begin{equation*}
\im (\varphi \overline{e^{i\sigma}})= d z_{\mu, \xi}.
\end{equation*}
This completes the proof.
\end{proof}

\section{The functional framework}

In the variational framework of the problem, solutions to
$(\ref{tp})$ can be found as critical points of the energy
functional $f_{\varepsilon}:E \rightarrow \R$ defined by
\begin{equation}
f_{\varepsilon}(u)=\frac{1}{2}\int_{\R^N} \left| \left(\frac
{\nabla}{i} - \varepsilon A(x)\right)u \right|^{2} dx +
\frac{\varepsilon^\alpha}{2}\int_{\R^N} V(x)|u|^2 dx
-\frac{1}{2^{*}}\int_{\R^N} |u|^{2^{*}} dx,
\end{equation}
on the real Hilbert space
\begin{equation}
E = D^{1,2} (\R^N,\C) = \left\{v \in L^{2^*}(\R^N,\C) \mid \int_{\R^N}
|\nabla v|^{2} dx < \infty \right\}
\end{equation}
endowed with the inner product
\begin{equation}
{\left\langle u, v \right\rangle }_{E} = \re \int_{\R^N} \nabla u \cdot
\overline {\nabla v} \, dx.
\end{equation}
We shall assume throughout the paper that
\begin{description}
\item[(N)] $N > 4$,
\item[(A1)] $A\in C^1(\R^N,\R^N) \cap L^{\infty}(\R^N,\R^N)\cap L^{r}(\R^N,\R^N)$ with
$1 < r < N$
\item[(A2)] $\operatorname{div}A \in  L^{N/2}(\R^N,\R)$,
\item[(V)] $V \in  C(\R^N,\R^N) \cap  L^\infty(\R^N, \R)\cap L^{s}(\R^N,\R)$, with $ 1 <
s < N/2$.
\end{description}
The functional $f_{\varepsilon}$ is well defined on $E$. Indeed,
\begin{equation*}
\int_{\R^N} \left| \left(\frac {\nabla }{i} - \varepsilon A(x)\right)u
\right|^{2}
= \int_{\R^N} | \nabla u|^{2} + \varepsilon^{2} \int_{\R^N} |A|^{2}
|u|^{2} - \re \int_{\R^N} \frac {\nabla u}{i} \cdot \varepsilon A\,
\overline{u},
\end{equation*}
and all the integrals are finite by virtue of \textbf{(A1)}.
Moreover, $f_\varepsilon \in C^{2}(E,\R)$.

In this section, we perform a finite--dimensional reduction on
$f_\varepsilon$ according to the methods of \cite{ab,am05}. Roughly
speaking, since the unperturbed problem (i.e. $(\ref{tp})$ with
$\ge=0$) has a whole $C^2$ manifold of critical points, we can
deform this manifold is a suitable manner and get a
\emph{finite--dimensional
  natural constraint} for the Euler--Lagrange functional associated to
$(\ref{tp})$. As a consequence, we can find solutions to
$(\ref{tp})$ in correspondence to (stable) critical points of an
auxiliary map
--- called the Melnikov function --- in finite dimension.

Now we focus on the case $\alpha =2$, as in the other cases $\alpha
\in [1,2[$ the magnetic potential $A$ no longer affects the
finite-dimensional reduction (see Remark $(\ref{fin})$).

So that we can write the functional $f_\varepsilon$ as
\begin{equation}
f_{\varepsilon}(u)= f_{0}(u)+ \varepsilon G_{1}(u)+\varepsilon^{2}
G_{2}(u)
\end{equation}
where
\begin{equation}
f_{0}(u)=\frac{1}{2}\int_{\R^N} |\nabla u
|^{2}-\frac{1}{2^{*}}\int_{\R^N} |u|^{2^{*}},
\end{equation}
\begin{equation}
G_{1}(u)=- \re \frac{1}{i} \int_{\R^N} \nabla u \cdot A\,
\overline{u}, \quad G_{2}(u)=\frac{1}{2}\int_{\R^N} |A|^{2} |u|^{2}
+ \frac{1}{2}\int_{\R^N} V(x) |u|^{2}.
\end{equation}
We can now use the arguments of \cite{ab,am05} to build a natural
constraint for the functional $f_\varepsilon$.

\begin{theorem}
Given $R > 0$ and $B_{R}=\left\{u \in E: || u|| \leq R \right\}$,
there exist $\varepsilon_{0}$ and a smooth function
$w=w(z,\varepsilon)=w(e^{i \sigma} z_{\mu, \xi}, \varepsilon)=
w(\sigma,\mu,\xi,\varepsilon)$, $w(z, \varepsilon):M=Z \cap
B_{R}\times (\varepsilon_{0},\varepsilon_{0}) \rightarrow E$ such that
\begin{enumerate}
\item $w(z,0)=0$  for all $z \in Z \cap B_{R}$
\item $w(z, \varepsilon)$ is orthogonal to $T_{z}Z$,
for all $(z, \varepsilon) \in M$. Equivalently $w(z, \varepsilon)
\in (T_{z}Z)^{\bot}$ \item the manifold $Z_{\varepsilon} =\left\{z+
w(z, \varepsilon): (z, \varepsilon) \in M\right\}$ is a natural
constraint for $f'_{\varepsilon}$: if $u \in Z_{\varepsilon}$
and $f'_{\varepsilon}|_{Z_{\varepsilon}}= 0$, then
$f'_{\varepsilon}(u)=0$.
\end{enumerate}
\end {theorem}
For future reference let us recall that $w$ satisfies 2. above and $D
f_{\varepsilon}(z+w) \in T_{z}Z$, namely $f''_0 (z)[w]+ \varepsilon
G'_1(z)+ o(\varepsilon) \in T_z Z$. As a consequence, if $G'_1 (z)
\bot T_z Z$ (to be proved as Lemma~\ref{lem:1.2}), one finds
\begin{equation}\label{eq:1}
w( \varepsilon, z)= - \varepsilon L_z G'_1 (z)+ o(\varepsilon),
\end{equation}
where $L_z$ denotes the inverse of the restriction to $(T_z Z)^{\bot}$
of $f''_0 (z)$.

\begin{lemma} \label{lem:1.2}
$G_1 (z)=0$ for all $z\in Z$.
\end{lemma}
\begin{proof}
\begin{eqnarray}
  G_{1}(z) &=& - \re \int_{\R^N}\frac{\nabla z}{i}\cdot A(x)\, \overline{z} \, dx= \left[ z= e^{ i \sigma} z_{\mu, \xi}\right]   \nonumber  \\
  &=& - \re \int_{\R^N} e^{ i \sigma}\frac{\nabla z_{\mu, \xi}}{i}\cdot A(x)\,
  e^{- i \sigma}\overline{z_{\mu, \xi}} \, dx= \nonumber \\
  &=& - \re \int_{\R^N}\frac{\nabla {z_{\mu, \xi}}}{i}\cdot A(x)\,  z_{\mu, \xi} \, dx= 0. \nonumber
\end{eqnarray}
\end{proof}
Hence we cannot hope to apply directly the tools contained in
\cite{ab}, since the \textit{Melnikov function} would vanish identically.
However, following \cite{am99}, we can find a slightly implicit
Melnikov function whose stable critical points produce critical points
of $f_\varepsilon$.

\begin{lemma}
Let $\Gamma : Z \rightarrow \R$  be defined by setting
\begin{equation} \label{eq:13}
\Gamma (z)= G_2 (z) - \frac{1}{2}\left(L_z G'_1 (z), G'_1 (z)\right).
\end{equation}
Then we have
\begin{equation} \label{eq:13bis}
f_ \varepsilon (z+ w (\varepsilon, z))= f_0 (z) + \varepsilon^2 \Gamma (z) + o(\varepsilon^2).
\end{equation}
\end{lemma}
\begin{proof}
Since $G_1|_Z \equiv 0$, then $G'_1 (z) \in (T_z Z)^\bot$. Then one finds
\begin{eqnarray*}
f_ \varepsilon (z+ w (\varepsilon, z)) &=& f_0 (z+ w (\varepsilon,
z))+ \varepsilon G_1 (z+ w (\varepsilon, z))+ \varepsilon^2 G_2(z+ w
(\varepsilon, z)) \\ &=& f_0 (z) + \frac{1}{2} f''_0 (z) [w,w]+
\varepsilon G_1 (z) + \varepsilon G'_1 (z) [w] + \varepsilon^2 G_2 (z)
+ o(\varepsilon^2).
\end{eqnarray*}
Using (\ref{eq:1}) and Lemma \ref{lem:1.2} the lemma follows.
\end{proof}

\begin{remark}
We notice that $\Gamma=G_2(z)+\frac{1}{2}\left(G'_1(z),\phi \right)$,
where $z$ stands for $e^{i \sigma}z_{\mu, \xi}$ and
$\phi=\lim_{\varepsilon \rightarrow 0} \frac{w}{\varepsilon}$.
\end{remark}

\begin{remark}
  By the definition of $z\,\in\,Z$, it results: $\Gamma (z)=
  \Gamma(e^{i \sigma}z_{\mu, \xi})= \Gamma(\sigma, \mu, \xi) $.  In
  the sequel, we will write freely $\Gamma(\sigma, \mu, \xi) \equiv
  \Gamma( \mu, \xi)$ since $\Gamma$ is $\sigma$-invariant.  Indeed, it
  is easy to check that $G_2$ is $\sigma$-invariant. In fact, by the
  definition of $G_2(z)$ and $z=e^{i \sigma } z_{\mu, \xi}$, it
  results:
\[
G_2(\sigma, \mu, \xi)=G_2 (e^{i \sigma } z_{\mu, \xi})=
\frac{1}{2}\int |A(x)|^2 |z_{\mu, \xi}|^2 \, dx + \frac{1}{2}\int
|V(x)| |z_{\mu, \xi}|^2 \, dx \equiv G_2( \mu, \xi).
\]
It remains to prove that $\left\langle G'_1 (z), \phi \right\rangle$
is $\sigma$-invariant. We will show that $\phi= e^{i \sigma } \psi
(\mu,\xi)$ with $\psi (\mu,\xi)\,\in\, \C$ independent on $\sigma$
which immediately gives
\begin{eqnarray*}
  \left\langle G'_1 ( e^{i \sigma } z_{\mu, \xi}), \phi  \right\rangle &=& - \re \int \frac{1}{i} e^{i \sigma } \nabla z_{\mu, \xi} \cdot A(x) e^{-i \sigma}\overline{\psi (\mu,\xi)}\, dx \\
  && - \re \int \frac{1}{i} \nabla \psi_{\mu, \xi} \cdot A(x)
  z_{\mu,\xi}\, dx = \left\langle G'_1 (z_{\mu, \xi}), \psi(\mu, \xi)
  \right\rangle.
\end{eqnarray*}
We begin to recall that $\phi= \lim_{\ge \rightarrow 0^+} \frac{w(\ge,
  z)}{\ge}$, where $w(\ge, z)$ is such that
\[
f'_\ge (e^{i \sigma } z_{\mu, \xi} + w(\sigma, \mu, \xi)) \,\in\,
T_{e^{i \sigma} z_{\mu,\xi}} Z.
\]
By (\ref{tang}), this condition means that
\begin{equation}\label{ipppp}
  f'_\ge (e^{i \sigma } z_{\mu, \xi} + w(\sigma, \mu,
  \xi))=\sum_{i=1}^{N} a_i e^{i \sigma} \frac{\partial
    z_{\mu,\xi}}{\partial \xi_i} + b e^{i \sigma} \frac{\partial
    z_{\mu,\xi}}{\partial \mu}+ d e^{i \sigma} i z_{\mu,\xi},
\end{equation}
with $a_1,\ldots,a_N$, $b$, $d, \,\in\,\R$.

Let $w(\sigma, \mu, \xi)= e^{i \sigma} \widetilde{w}$ with
$\widetilde{w}\, \in \, D^{1,2}(\R^N, \C)$. Testing $(\ref{ipppp})$ by
$e^{i \sigma} v(x)$ with $v \, \in \, D^{1,2}(\R^N, \C)$, we derive
that $z_{\mu,\xi} + \widetilde{w}$ is a solution of an equation
independently on $\sigma$. Thus, also $\widetilde{w}$ is independent
on $\sigma$ and it can be denoted as $\widetilde{w}(\mu, \xi)$. Set
$\psi( \mu, \xi)= \lim_{\ge \rightarrow 0^+} \frac{ \widetilde{w}(
  \mu, \xi)}{\ge} $, we deduce that $\phi= e^{i \sigma } \psi
(\mu,\xi)$.
\end{remark}

\section{Asymptotic study of \(\Gamma\)}

In order to find critical points of $\Gamma$ it is convenient to study
the behavior of $\Gamma$ as $ \mu \rightarrow 0$ and as $\mu + |\xi|
\rightarrow \infty$.
Our goal is to show:

\begin{proposition} \label{prop:1}
  $\Gamma$ can be extended smoothly to the hyperplane $\left\{(0,
    \xi) \in \R \times \R^N \right\}$ by setting
\begin{equation}
\Gamma (0, \xi)=0.
\end{equation}
Moreover there results
\begin{equation}
  \Gamma (\mu, \xi) \rightarrow 0, \quad\hbox{as $\mu +|\xi| \rightarrow +
\infty$}.
\end{equation}
\end{proposition}
The proof of this Proposition is rather technical, so we split it into
several lemmas in which we will use the formulation of
$\Gamma=G_2(z)+\frac{1}{2}\left(G'_1(z),\phi \right)$,
where
$\phi=\lim_{\varepsilon \rightarrow 0} \frac{w}{\varepsilon}$.

\begin{lemma}
Under assumption \textbf{(A1)} there holds
\begin{equation} \label{eq:5}
\lim_{\mu \rightarrow 0^+} \frac{1}{2}\int_{\R^N}|A(x)|^2
| z_{\mu ,\xi} |^2 dx = 0.
\end{equation}
\end{lemma}
\begin{proof}
Let $z=e^{i \sigma} z_{\mu ,\xi} \in Z$. Then
\begin{eqnarray}
H_2(z) &=&  \frac{1}{2}\int_{\R^N}|A(x)|^2 | z_{\mu ,\xi} |^2 dx \label{eq:2}
\\ &=&
\frac{1}{2}\int_{\R^N}|A(x)|^2 {\left( \kappa_N
\mu^{-\frac{(N-2)}{2}} \left(1+{\left|\frac{x-\xi}{\mu}
\right|}^2 \right)^\frac{2-N}{2}\right) }^2 dx \notag \\ &=&
\frac{\kappa_{N}^{2}}{2 \mu^{(N-2)}} \int_{\R^N}
\frac{|A(x)|^2}{\big(1+{\big|\frac{x-\xi}{\mu}
\big|}^{2} \big)^{N-2}} \, dx \notag
\end{eqnarray}
Using the change of variable $y=\frac{x-\xi}{\mu}$, or $x=\mu y+ \xi$,
we can write
\begin{eqnarray*}
H_2(z) &=& \frac{ \kappa_{N}^{2}}{2 \mu^{(N-2)}} \int_{\R^N}|A(\mu
y+\xi)|^2 \frac{1}{(1+ |y|^{2})^{N-2}}\mu^{N} dy \\ &=& \frac{
\kappa_{N}^{2}}{2} \mu^{2} \int_{\R^N}\frac{|A(\mu y+\xi)|^2}{(1+
|y|^{2})^{N-2}} dy
\end{eqnarray*}
and using the hypothesis \textbf{(A1)}
\begin{equation}
  H_2(z) \leq \mu^2 C_N {{\left\|A\right\|}^{2}_{\infty}} \int_{\R^N}
  \frac{1}{(1+ |y|^{2})^{N-2}} \, dy,
\end{equation}
the lemma follows.
\end{proof}

The proof of the following Lemma is similar and thus omitted.
\begin{lemma}
Under assumption $\textbf{(V)}$ there holds
\begin{equation} \label{eq:5bis}
\lim_{\mu \rightarrow 0^+} \frac{1}{2} \int_{\R^N}V(x) |
z_{\mu ,\xi} |^2 dx = 0.
\end{equation}
\end{lemma}

\begin{lemma}
There holds
\begin{equation} \label{eq:6}
\lim_{\mu \to 0^+} \langle G_1'(z),\phi \rangle =0.
\end{equation}
\end{lemma}
\begin{proof}
We write
\begin{equation*}
\left\langle G'_1(z), \phi \right\rangle_ E =\alpha_1+ \alpha_2,
\end{equation*}
where
\begin{eqnarray}
\alpha_1 &=&- \re \int_{\R^N}\frac{\nabla z}{i} \cdot A(x) \overline \phi \,
dx \\ \alpha_2 &=& - \re \int_{\R^N} \frac{\nabla \phi}{i} \cdot A(x)
\overline z \, dx.
\end{eqnarray}
It is convenient to introduce $\phi^*(y)$ by setting
\begin{equation*}
\phi^*(y)= \phi^{*}_{\mu, \xi}(y)= \mu^{\frac{N}{2}-1} \phi (\mu y +
\xi)
\end{equation*}
Using the expression of $z= e^{i \sigma} \mu^{-\frac{(N-2)}{2}}z_0
(\frac{x-\xi}{\mu})$ and the change of variable $x=\mu y + \xi$ we can
write:
\begin{eqnarray*}
\alpha_1 &=& - \re \int_{\R^N}\frac{1}{i}\nabla_{x} e^{i \sigma}
\mu^{-\frac{(N-2)}{2}}z_0 \left(\frac{x-\xi}{\mu}\right) \cdot A(x)
\overline {\phi (x)} \, dx \\ &=& - \re \int_{\R^N}\frac{1}{i} e^{i \sigma}
\nabla_{y} z_0 (y) \mu^{\frac{N}{2}} \cdot A(\mu y + \xi) \overline
{\phi (\mu y + \xi)} \, dy \\ &=& -\mu \re \int_{\R^N}\frac{1}{i} e^{i \sigma}
\nabla_y z_0 (y) \cdot A(\mu y + \xi) \overline {\phi^{*}(y)} \, dy
\end{eqnarray*}
and
\begin{eqnarray*}
\alpha_2 &=& - \re \int_{\R^N}\frac{1}{i}\nabla_{x} \phi (x) \cdot A(x)
e^{-i \sigma} \mu^{-\frac{(N-2)}{2}}z_0 \left(\frac{x-\xi}{\mu}\right)
\, dx \\ &=& - \re \int_{\R^N}\frac{1}{i}\nabla_{y} \phi (\mu y+ \xi)
\mu^{-1} \cdot A(\mu y+\xi) e^{-i \sigma} \mu^{-\frac{(N-2)}{2}} \mu^N
z_0 (y) \, dy \\ &=& - \re \int_{\R^N}\frac{1}{i}\nabla \phi (\mu y+ \xi)
\mu^{-1} \cdot A(\mu y+\xi) e^{-i \sigma} \mu^{\frac{N}{2}+1} z_0 (y)
\, dy \\ &=& -\mu \re \int_{\R^N}\frac{1}{i}\nabla \phi^{*} (y) \cdot A(\mu
y+\xi) e^{-i \sigma} z_0 (y) \, dy.
\end{eqnarray*}
Now the conclusion follows easily from the next lemma.
\end{proof}

\begin{lemma}\label{tom}
As $\mu \to 0^{+}$,
\begin{equation} \label{eq:4}
\phi^{*}_{\mu,\xi} \rightarrow 0 \quad\text{strongly in $E$}.
\end{equation}
\end{lemma}

\begin{proof}
For all $v \in E$, due to the divergence theorem, we have
\begin{eqnarray*}
\left\langle G'_1(z), v \right\rangle_ E &=& - \re \int_{\R^N}\frac{\nabla z}{i}
\cdot A(x) \overline v \, dx - \re \int_{\R^N} \frac{\nabla v}{i} \cdot A(x)
\overline z \, dx\\
&=& - \re \int_{\R^N}\frac{\nabla z}{i} \cdot A(x)
\overline v \, dx - \re \int_{\R^N}\frac{1}{i} \sum_{j=1}^{N}\frac{\partial
v}{\partial x_{j}}A_{j} (x) \overline{z} \, dx \\
&=& - \re \int_{\R^N}
\frac{\nabla z}{i} \cdot A(x) \overline v \, dx + \re \int_{\R^N}\frac{1}{i}
\sum_{j=1}^{N} v \frac{\partial}{\partial x_{j}}(A_{j} \overline{z})
\, dx\\
&=&- \re \int_{\R^N} \frac{\nabla z}{i} \cdot A(x) \overline v \, dx +
\re \int_{\R^N}\frac{1}{i} v \operatorname{div} A \, \overline{z} \, dx + \re
\int_{\R^N}\frac{1}{i} v A \cdot \overline{\nabla z} \, dx \\
&=& - 2 \re \int_{\R^N}
\frac{\nabla z}{i} \cdot A(x) \overline v \, dx - \re \int_{\R^N}\frac{1}{i}
 \operatorname{div} A \,z \overline{v} \, dx
\end{eqnarray*}
where the last integral is finite by assumption \textbf{(A2)} and
\begin{multline}
(f''_{0} (z) w_{\mu,\xi},v)=
\re \int_{\R^N}\nabla  w_{\mu,\xi} \cdot \, \overline{\nabla v} \, dx- \re
\int_{\R^N}|z|^{2^{*}-2}  w_{\mu,\xi} \overline v \, dx  \\
- \re \int_{\R^N}(2^{*}-2)|z|^{2^{*}-4} \re (z \overline w_{\mu,\xi} )z
\overline v \, dx .
\end{multline}
We know that $w_{\mu,\xi}=-\varepsilon L_{ e^{i \sigma}  z_{\mu,\xi}} G'_1
(  e^{i \sigma}     z_{\mu,\xi})+ o(\varepsilon)$, and hence
\begin{equation}
\langle f''_{0}(z) \phi_{\mu,\xi}, v \rangle_{E} = -\langle G'_1
(z),v \rangle_{E}, \quad \forall v \in E
\end{equation}
where $\phi_{\mu,\xi} = \lim_{\epsilon \to 0}
\frac{w_{\mu,\xi}}{\epsilon}$. This implies that $\phi_{\mu,\xi}$
solves
\begin{multline*}
\re \int_{\R^N}\nabla  \phi_{\mu,\xi} \cdot \, \overline{\nabla v} dx
- \re \int_{\R^N}|z|^{2^{*}-2}  \phi_{\mu,\xi} \overline v \, dx -
\re \int_{\R^N}(2^{*}-2)|z|^{2^{*}-4} \re (z \overline
\phi_{\mu,\xi} )z \overline v \, dx  \\
=2  \re\int_{\R^N}\frac{1}{i} \nabla z \cdot A(x) \overline {v} dx + \re
\int_{\R^N}\frac{1}{i}  \operatorname{div} A  \,z \overline{v} \, dx.
\end{multline*}
Multiplying by $\mu^{\frac{N}{2}-1}$ and using the expression of $z= e^{i
  \sigma} \mu^{-\frac{(N-2)}{2}} z_0 (\frac{x-\xi}{\mu})$, we get
\begin{multline*}
  \re \int_{\R^N}\mu^{\frac{N}{2}-1} \nabla_{x} \phi_{\mu,\xi}(x)
  \overline{\nabla v}dx -\re \int_{\R^N}\mu^{-2} \left|\,z_0
    \left(\frac{x-\xi}{\mu}\right)\,\right|^{2^*
    -2} \mu^{\frac{N}{2}-1} \phi_{\mu,\xi}(x)\overline{v} dx  \\
  -\re \int_{\R^N}(2^* -2) \mu^{N-4}\left|\,z_0
    \left(\frac{x-\xi}{\mu}\right)\,\right|^{2^* -4} \re \left( e^{i
      \sigma}\mu^{-\frac{N}{2}+1}z_{0} \left(\frac{x-\xi}{\mu}
    \right) \mu^{\frac{N}{2}-1} \overline{\phi_{\mu,\xi}(x)} \right)\times
  \\
  \times e^{i \sigma}\mu^{-\frac{N}{2}+1}
  z_{0} \left(\frac{x-\xi}{\mu}\right)\overline{v}dx \\
  =2 \re \int_{\R^N}\frac{1}{i} e^{i \sigma} \nabla_{x} z_{0} \left(\frac{x-\xi}{\mu}
  \right) \cdot A(x) \overline {v} dx + \re \int_{\R^N}\frac{1}{i}  \operatorname{div} A \, e^{i \sigma}z_{0} \left(\frac{x-\xi}{\mu} \right)\overline {v} dx .
\end{multline*}
Using the expression of
$\phi^{*}(\frac{x-\xi}{\mu})=\mu^{\frac{N}{2}-1} \phi_{\mu,\xi}(x)$,
we have
\begin{multline*}
  \re \int \nabla_{x} \phi^{*} \left(\frac{x-\xi}{\mu}\right)
  \overline{\nabla v}dx -\re \int \mu^{-2} \left|\,z_0
    \left(\frac{x-\xi}{\mu}\right)\,\right|^{2^* -2} \phi^{*}
  \left(\frac{x-\xi}{\mu}\right)
  \overline{v} dx \\
  - \re \int (2^* -2) \mu^{N-4}\left|\,z_0
    \left(\frac{x-\xi}{\mu}\right)\,\right|^{2^* -4}
  \re \left( e^{i \sigma}\mu^{-\frac{N}{2}+1}z_{0} \left(\frac{x-\xi}{\mu}
  \right)  \overline{\phi^{*} \left(\frac{x-\xi}{\mu}\right)}
  \right)  \times  \\
  \times e^{i \sigma} \mu^{-\frac{N}{2}+1}
  z_{0} \left(\frac{x-\xi}{\mu}\right) \overline{v} dx \\
  =2 \re \int \frac{1}{i} e^{i \sigma} \nabla_{x} z_{0} \left(\frac{x-\xi}{\mu}
  \right) \cdot A(x) \overline {v} dx + \re \int \frac{1}{i}  \operatorname{div} A(x) \,e^{-i \sigma}z_{0} \left(\frac{x-\xi}{\mu} \right)  \overline{v}dx .
\end{multline*}
then, the change of variable $x=\mu y+\xi$ yields
\begin{multline*}
  \re \int \mu^{-2} \nabla_{y} \phi^{*}(y) \overline{\nabla_{y} v(\mu
    y+\xi)} \mu^{N} dy -\re \int \mu^{N-2}
  \left|\,z_0(y)\,\right|^{2^* -2} \phi^{*}(y)
  \overline{v(\mu y+\xi)} dy  \\
  {} - \re \int (2^* -2) \mu^{N-4}\left|\,z_0(y) \right|^{2^* -4} \re
  \left( e^{i \sigma}\mu^{2(-\frac{N}{2}+1)}z_{0}(y)
    \overline{\phi^{*}(y)} \right) e^{i \sigma}
  z_{0}(y) \overline{v(\mu y+\xi)} \mu^{N}dy \\
  = 2 \re \int \frac{1}{i} e^{i \sigma} \nabla_{y} z_{0}(y) \cdot
  A(\mu
  y+\xi)) \overline {v(\mu y+\xi)} \mu^{N-1}dy \\
  \qquad \qquad + \re \int \frac1i  \operatorname{div}
  A(\mu y+\xi) \, e^{-i \sigma}z_0(y) \overline{v(\mu y+\xi)} \mu^{N}dy .
\end{multline*}
Replacing $x=y$ and dividing by $\mu^{N-2}$, it results
\begin{multline*}
  \re \int_{\R^N} \nabla_{x} \phi^{*}(x) \overline{\nabla_{x} v(\mu
    x+\xi)} \, dx - \re \int_{\R^N} \left|\,z_0(x)\,\right|^{2^* -2}
  \phi^{*}(x)\overline{v(\mu x+\xi)} \, dx \\
  \qquad {} - \re \int_{\R^N}(2^* -2) \left|\,z_0(x) \right|^{2^* -4}
  \re \left( e^{i \sigma}z_{0}(x) \phi^{*}(x) \right) e^{i \sigma}
  z_{0}(x) \overline{v(\mu x+\xi)} \, dx \\
  =2 \mu \re \int_{\R^N}\frac{1}{i} e^{i \sigma} \nabla_{x} z_{0}(x)
  \cdot A(\mu
  x+\xi)) \overline {v(\mu x+\xi)} \, dx  \\
  \qquad {} + \mu^2 \re \int_{\R^N}\frac{1}{i}
  \operatorname{div}_{y} A(\mu x+\xi) \, e^{i \sigma} {z_{0}(x)} \overline{ v(\mu x+\xi)}
  \, dx.
\end{multline*}
This means that, if we write $\tau_{\mu,\xi}(x) = \mu x + \xi$,
\[
\left\langle f''_0 (e^{i \sigma}z_0) \phi^* , v \circ \tau_{\mu, \xi} \right\rangle= \int k_{\mu,\xi} \overline{ v \circ \tau_{\mu, \xi}  }
\]
for all test function $v$, in particular that
\[
 f''_0 (e^{i \sigma} z_0) \phi^*=  k_{\mu,\xi}
\]
where
\[
k_{\mu, \xi}(x) = \frac{2}{i} \mu e^{i \sigma} \nabla_{x} z_{0}(x)
\cdot A(\mu x+\xi) + \frac{1}{i} \mu^2 e^{i \sigma}
\operatorname{div}_y A(\mu x+\xi) \, z_{0}(x).
\]
We conclude that $\phi^{*}$ is a solution of
\begin{equation}
  \phi^{*}(x)= L_{ e^{i \sigma} z_{0}} k_{\mu, \xi}(x)
\end{equation}
Our assumptions on $A$ (i.e. \textbf{(A1)} and \textbf{(A2)}) imply
immediately that
\begin{equation}
  k_{\mu, \xi} \rightarrow 0 \quad \hbox{in $E$ as $\mu \rightarrow 0$}.
\end{equation}
From the continuity of $L_{e^{i \sigma}z_0}$ we deduce that
\begin{equation}
  \lim_{\mu \to 0+} \phi^* =  \lim_{\mu \to 0+} L_{e^{i \sigma}z_0}k_{\mu,\xi}
= 0.
\end{equation}
This completes the proof of the Lemma.
\end{proof}

\begin{lemma}
  Under assumption \textbf{(A1)}, there holds
\[
\lim_{\mu + |\xi| \to + \infty } H_2 (\mu,\xi) = 0,
\]
where $H_2$ is defined in (\ref{eq:2}).
\begin{proof}
  Firstly, assume that $\mu \to \overline{\mu} \in (0, +\infty)$ and
  $\mu + |\xi| \to + \infty$. We notice that
\begin{eqnarray*}
  H_2 (\mu,\xi) &=& \frac{\mu^{-(N-2)}}{2} \int_{\R^N} |A(x)|^2 z_0^2 \left(\frac{x-\xi}{\mu}       \right) dx\\
  &=& \frac{\mu^{-(N-2)}}{2}\int_{|x| \leq \frac{|\xi|}{2}} |A(x)|^2 z_0^2 \left(\frac{x-\xi}{\mu} \right) dx \\
  && {} +\frac{\mu^{-(N-2)}}{2}\int_{|x| > \frac{|\xi|}{2}} |A(x)|^2 z_0^2 \left(\frac{x-\xi}{\mu} \right) dx.
\end{eqnarray*}
Moreover,
\begin{eqnarray*}
  &&\frac{\mu^{-(N-2)}}{2} \int_{|x| \leq \frac{|\xi|}{2}} |A(x)|^2 z_0^2 \left(\frac{x-\xi}{\mu} \right) dx \\
  && \qquad \leq \frac{\mu^{-(N-2)}}{2}||A||^2_{\infty} \omega_{N} \frac{|\xi|^N}{2^N} \sup_{|x| \leq \frac{|\xi|}{2}} z_0^2 \left(\frac{x-\xi}{\mu} \right) \\
  && \qquad = \frac{\mu^{-(N-2)}}{2}||A||^2_{\infty} \omega_{N} \frac{|\xi|^N}{2^N} \sup_{|x| \leq \frac{|\xi|}{2}} \frac{k_N^2 \mu^{2(N-2)}}{\left[ \mu^2 + | x - \xi|^2\right]^{N-2}} \\
  && \qquad \leq \frac{\mu^{-(N-2)}}{2}||A||^2_{\infty} \omega_{N} \frac{|\xi|^N}{2^N} \sup_{|x| \leq \frac{|\xi|}{2}} \frac{k_N^2}{\left[ \mu^2 + \left|\, |x| - |\xi|\, \right|^2\right]^{N-2}}  \\
  && \qquad \leq \frac{\mu^{-(N-2)}}{2}||A||^2_{\infty} \omega_{N} \frac{|\xi|^N}{2^N} \frac{k_N^2}{\left[ \mu^{2} + \frac{|\xi|^2}{4}\right]^{N-2}},
\end{eqnarray*}
where $\omega_N$ is the measure of $S^{N-1}= \left\{x \in \R^N : |x|
=1 \right\}$.  Since $N > 4$,  we infer
\begin{equation*}
  \frac{k_N^2 |\xi|^N}{\left[ \mu^{2} + \frac{|\xi|^2}{4}\right]^{N-2}} \to 0 \quad \hbox{as $|\xi| \to +\infty$.}
\end{equation*}
Finally, we deduce
\begin{equation*}
  \frac{\mu^{-(N-2)}}{2} \int_{|x| \leq \frac{|\xi|}{2}} |A(x)|^2 z_0^2 \left(\frac{x-\xi}{\mu} \right) dx \to 0
\end{equation*}
as $\mu \to \overline{\mu}$ and $|\xi| \to +\infty$.

On the other hand, we have
\begin{eqnarray*}
&&\frac{\mu^{-(N-2)}}{2} \int_{|x| > \frac{|\xi|}{2}} |A(x)|^2 z_0^2 \left(\frac{x-\xi}{\mu} \right) dx \\
&& \qquad\leq \frac{\mu^{-(N-2)}}{2} ||A||_{\infty}^{2}\int_{|x| > \frac{|\xi|}{2}}  z_0^2 \left(\frac{x-\xi}{\mu} \right) dx \\
&& \qquad = \frac{\mu^{N-(N-2)}}{2} ||A||_{\infty}^{2}\int_{|\mu x +\xi| > \frac{|\xi|}{2}}  z_0^2(x) dx.
\end{eqnarray*}
Since $z_0^2 \in L^1 (\R^N)$, we deduce that
\begin{equation*}
\frac{\mu^{2}}{2} ||A||_{\infty}^{2}\int_{|\mu x +\xi| > \frac{|\xi|}{2}}  z_0^2(x) dx \to 0
\end{equation*}
as  $\mu \to \overline{\mu}$ and $|\xi| \to +\infty$, and thus
\begin{equation*}
\frac{\mu^{-(N-2)}}{2} \int_{|x| > \frac{|\xi|}{2}} |A(x)|^2 z_0^2 \left(\frac{x-\xi}{\mu} \right) dx \to 0
\end{equation*}
as  $\mu \to \overline{\mu}$ and $|\xi| \to +\infty$.\\
Finally, we can conclude that $H_2 (\mu, \xi) \to 0$ as $\mu \to \overline{\mu}$ and $|\xi| \to +\infty$.\\
Conversely, assume that $\mu \to +\infty$. After a suitable change of variable, it results
\begin{equation*}
H_2(\mu,\xi)= \frac{\mu^2}{2} \int_{\R^N} |A(\mu y + \xi)|^2
|z_0(y)|^2 dy.
\end{equation*}
By assumption \textbf{(A1)}, we can fix $1 < r < \frac{N}{2}$ such that $A^2 \in L^r (\R^N)$. Moreover, let be $s= \frac{r}{r-1}$. It is immediate to check that $2s > 2^*$ and then $z_0^{2s} \in L^1 (\R^N)$. By (A1) and Holder inequality, we deduce that
\begin{multline*}
 \int_{\R^N} |A(\mu y + \xi)|^2 |z_0(y)|^2 dy \\
 \leq  \left(\int_{\R^N} |A(\mu y + \xi)|^{2r}  dy \right)^{\frac{1}{r}}  \left(\int_{\R^N} |z_0(y)|^{2s} dy \right)^{\frac{1}{s}} \\
 \leq \mu^{-\frac{N}{r}}\left(\int_{\R^N} |A( y)|^{2r}  dy \right)^{\frac{1}{r}} \left(\int_{\R^N}|z_0(y)|^{2s} dy \right)^{\frac{1}{s}}.
\end{multline*}
As a consequence, by the above inequality, we infer for $\mu$ small
\begin{multline*}
G_2(\mu,\xi)= \frac{\mu^2}{2} \int_{\R^N} |A(\mu y + \xi)|^2 |z_0(y)|^2 dy \\
\leq \mu^{2-\frac{N}{r}}\left(\int_{\R^N} |A( y)|^{2r}  dy \right)^{\frac{1}{r}} \left(\int_{\R^N}|z_0(y)|^{2s} dy \right)^{\frac{1}{s}}.
\end{multline*}
Now, we notice that $ r < \frac{N}{2}$ implies $2- \frac{N}{r} < 0$
and thus by the above inequality we can conclude that $G_2(\mu, \xi)$
tends to $0$ as $\mu \to +\infty$.
\end{proof}
\end{lemma}

Arguing as before we can deduce the following result.

\begin{lemma}
Under assumption \textbf{(V)}, there holds
\[
\lim_{\mu + |\xi| \to + \infty } \int_{\R^N} V(x)
|z_{\mu,\xi}(x)|^2 \, dx = 0.
\]
\end{lemma}
In order to describe the behavior of the term $\left\langle G'_1(z),
\phi \right\rangle_ E$ as $\mu + |\xi| \to + \infty$, we need the
following lemma.

\begin{lemma} \label{lemma:4.8}
 There is a constant $C_N >0$ such that
 \begin{equation}
{ \left\| \, \phi \, \right\|}_E \, \leq \, C_N \, \quad\text{ for
all $\mu>0$ and  for all $\xi \in \R^N$}.
 \end{equation}
\end{lemma}
 \begin{proof}
 We know that for all $\ge > 0$ and all $z \in Z$
 \begin{equation*}
 w(\ge, z)=- L_z G'_1(z) + o(\ge)
 \end{equation*}
 so that
 \begin{equation*}
 \phi = \lim_{ \ge \rightarrow 0} \frac{ w(\ge, z)}{\ge}= -L_z G'_1(z)
 \end{equation*}
 and
  \begin{equation*}
    {\left\|\phi \right\|}_E \, \leq \, \left\| L_z \right\|\,\left\| G'_1(z)\right\|.
 \end{equation*}
 We claim that $\|L_z\|$ is bounded above by a constant independent of
 $\mu$ and $\xi$.  Indeed:
 \begin{align*}
   & \left\| L_z \right\|= \sup_{  \left\| \varphi \right\|=1 }   \left\| L_z \varphi \right\|=\sup_{  \substack{\left\| \varphi \right\| =1 \\ \left\| \psi \right\| =1}}   \left| \left\langle L_z \varphi, \psi \right\rangle  \right|\\
   & = \sup_{ \substack{\left\| \varphi \right\| =1 \\ \left\| \psi
       \right\| =1} } \left| \int_{\R^N}\nabla \varphi \cdot
     \overline{\nabla \psi} - \re
     \int_{\R^N}\left|z_{\mu,\xi}\right|^{2^*-2} \varphi \overline{\psi} \right. \\
& \hspace{3cm}
     \left. {} - (2^*-2) \int_{\R^N}\left|z_{\mu,\xi}\right|^{2^*-4}\re (\overline{\varphi} z_{\mu,\xi})   \re (\overline{\psi} z_{\mu,\xi}) \right| \\
   & \leq\, \sup_{ \substack{\left\| \varphi \right\| =1 \\ \left\| \psi
       \right\| =1} } \left( \int_{\R^N}\left|\nabla \varphi \right| \,
     \left|\overline{\nabla \psi}\right| + \re
     \int_{\R^N}\left|z_{\mu,\xi}\right|^{2^*-2} \left|\varphi \right| \,
     \left|\overline{\psi}\right|
     + (2^*-2) \int_{\R^N}\left|z_{\mu,\xi}\right|^{2^*-2}\left|\varphi\right| \, \left|\overline{\psi}\right| \right) \\
   & \leq\,
   \sup_{  \substack{\left\| \varphi \right\| =1 \\ \left\| \psi \right\| =1} } \left( \int_{\R^N}\left|\nabla \varphi \right| \, \left|\overline{\nabla \psi}\right| + (2^*-1) \int_{\R^N}\left|z_{\mu,\xi}\right|^{2^*-2}\left|\varphi\right| \, \left|\overline{\psi}\right| \right) \\
   & \leq\, \sup_{ \substack{\left\| \varphi \right\| =1 \\ \left\| \psi
       \right\| =1} } \left( \int_{\R^N}\left|\nabla \varphi \right|^2
   \right)^{1/2} \left( \int_{\R^N}\left|\nabla \psi \right|^2
   \right)^{1/2}+ (2^*-1)
   \left(\int_{\R^N}\left|z_{\mu,\xi}\right|^{2^*}\right)^{(2^*-2)/2^*}\times
   \\
   &\qquad \qquad\qquad \times
   \left(\int_{\R^N}\left|\varphi\right|^{2^*}\right)^{1/2^*}\left(\int_{\R^N}\left|\psi\right|^{2^*}\right)^{1/2^*}
  \end{align*}
We observe that
 \begin{eqnarray*}
   \left(\int_{\R^N}\left|z_{\mu,\xi}\right|^{2^*}\right)^{1/2^*} &=& \mu^{-\frac{(N-2)}{2}}
   \left(\int_{\R^N}\left|z_0 \left(\frac{x-\xi}{\mu} \right)\right|^{2^*}\right)^{1/2^*}\\
   &=& \mu^{-\frac{(N-2)}{2}} \left(\int_{\R^N}\left|z_0(y)\right|^{2^*} \mu^N \, \right)^{1/2^*}={\left\|    z_0\right\|}_{L^{2^*}}.
 \end{eqnarray*}
Hence
 \begin{eqnarray*}
   \left\| L_z \right\|&\leq&  \sup_{  \substack{\left\| \varphi \right\| =1 \\ \left\| \psi \right\| =1} } \left(1+ (2^*-1) {\left\|\, z_0 \right\|\,}^{(2^*-2)}_{L^{2^*}}{\left\| \varphi \right\|}_{L^{2^*}}
     {\left\| \psi \right\|}_{L^{2^*}}\right)\\
   &\leq&   \sup_{  \substack{\left\| \varphi \right\| =1 \\ \left\| \psi \right\| =1} } \left(1+ (2^*-1)C'_N {\left\| z_0 \right\|}^{(2^*-2)}_{E}{\left\| \varphi \right\|}_{E}
     {\left\| \psi \right\|}_{E}\right)\\
   &\leq&  1+ (2^*-1)C'_N {\left\| z_0 \right\|}^{(2^*-2)}_{E} \equiv C^1_N
 \end{eqnarray*}
 where $ C^1_N$ is a constant independent from $\mu$ and $\xi$. At
 this point it results:
\begin{equation*}
  \left\|\phi \right\|\leq C^1_N \left\|G'_1(z) \right\|
\end{equation*}
and we have to evaluate $\left\|G'_1(z) \right\|$ :
\begin{eqnarray*}
\left\|G'_1(z) \right\|&=& \sup_{\left\|\varphi\right\|=1} \left|   \left\langle G'_1(z), \varphi  \right\rangle \right|\\
&=& \sup_{\left\|\varphi\right\|=1} \left| \left( -\re
\int_{\R^N}\frac{\nabla z}{i} \cdot A(x) \overline{\varphi}\, dx -\re
\int_{\R^N}\frac{\nabla \phi}{i} \cdot A(x) \overline{z}\, dx
\right)
\right|\\
&\leq& \sup_{\left\|\varphi\right\|=1}
\left(  \int_{\R^N}\left|\nabla z_{\mu,\xi}\right| \left| A(x)\right| \left|\overline{\varphi}\right|\, dx + \int_{\R^N} \left| \nabla \varphi\right|  \left|A(x)\right| \left| z_{\mu,\xi} \right|\, dx    \right)\\
&\leq&{ \left\|A \right\|}_{L^N}
\sup_{\left\|\varphi\right\|=1}\left( { \left\|z_0 \right\|}_E \,{ \left\| \varphi \right\|}_E  C''_N     \right) \\
&\leq&  { \left\|A \right\|}_{L^N} { \left\|z_0 \right\|}_E  C''_N
\equiv C^2_N
\end{eqnarray*}
with $C^2_N$ independent from $\mu$ and $\xi$.\\
Finally,
\begin{equation*}
\left\|\, \phi \, \right\|\,\leq\,C^1_N  C^2_N \equiv C_N
\end{equation*}
with $C_N$ independent from $\mu$ and $\xi$ and the lemma is proved.
\end{proof}

\begin{remark}
It is easy to check that $\left\|\phi^* \right\|=\left\|\phi \right\|$.
\end{remark}

\begin{lemma}
There holds
\[
\lim_{\mu + |\xi| \to + \infty } \left\langle G'_1(z), \phi
\right\rangle_ E = 0.
\]
\end{lemma}
\begin{proof}
Firstly, assume that $\mu \to \overline{\mu} \in (0, +\infty)$ and
$\mu + |\xi| \to + \infty$. We can write
\[ \left\langle G'_1(z), \phi \right\rangle_ E=\alpha_1+\alpha_2
\]
where
\begin{eqnarray}
\alpha_1 &=&- \re \int_{\R^N}\frac{\nabla z}{i} \cdot A(x) \overline \phi \,
dx \\
\alpha_2 &=& - \re \int_{\R^N} \frac{\nabla \phi}{i} \cdot A(x)
\overline z \, dx.
\end{eqnarray}
Using the expression of $z= e^{i \sigma} \mu^{-\frac{(N-2)}{2}}z_0
(\frac{x-\xi}{\mu})$ and by assumption \textbf{(A1)} and the H\"older
inequality we have:
\begin{eqnarray*}
\alpha_1 &=& - \re \int_{\R^N}\frac{1}{i}\nabla_{x} e^{i \sigma}
\mu^{-\frac{(N-2)}{2}}z_0 \left(\frac{x-\xi}{\mu}\right) \cdot A(x) \overline
\phi (x) \, dx \\
&\leq & \mu^{-\frac{(N-2)}{2}} \left( \int_{\R^N}\left| \nabla_{x}z_0 \left(\frac{x-\xi}{\mu}\right)\right|^2 \, dx \right)^{1/2} \left(\int_{\R^N}\left( |A(x)|\, | \overline{\phi}| \right)^2 \, dx   \right)^{1/2}\\
& \leq & \mu^{-\frac{(N-2)}{2}}\left\| A \right\|_{L^N (\R^N)}\, \left\|\phi \right\|_{L^{2^*} (\R^N)}\, \left( \int_{\R^N}\left| \nabla_{x}z_0 \left(\frac{x-\xi}{\mu}\right)\right|^2 \, dx \right)^{1/2}
\end{eqnarray*}
We notice that
\begin{multline*}
\int_{\R^N} \left| \nabla_{x}z_0 \left(\frac{x-\xi}{\mu}\right)\right|^2 \, dx = \int_{|x|\, \leq \, |\xi|/2} \left| \nabla_{x}z_0 \left(\frac{x-\xi}{\mu}\right)\right|^2 \, dx \\
{} +  \int_{|x|\, > \, |\xi|/2} \left| \nabla_{x}z_0 \left(\frac{x-\xi}{\mu}\right)\right|^2 \, dx
\end{multline*}
and
\begin{equation*}
\left| \nabla_{x}z_0 \left(\tfrac{x-\xi}{\mu}\right)\right|^2= \mu^{2(2-N)} (2-N)^2 \kappa^2_N \frac{|x-\xi|^2}{\left( \mu^2 + |x-\xi|^2 \right)^N}.
\end{equation*}
Moreover, setting $C^2_N:=(2-N)^2 \kappa^2_N$,
\begin{eqnarray*}
\int_{|x|\, \leq \, |\xi|/2} \left| \nabla_{x}z_0 \left(\frac{x-\xi}{\mu}\right)\right|^2 \, dx& \leq & \omega_N \frac{|\xi|^N}{2^N} \sup_{  |x|\, \leq \, |\xi|/2  } \, \left| \nabla_{x}z_0 \left(\frac{x-\xi}{\mu}\right)\right|^2 \\
  &=& \omega_N \frac{|\xi|^N}{2^N} \sup_{|x|\, \leq \, |\xi|/2}\mu^{2(2-N)} (2-N)^2 \kappa^2_N \frac{|x-\xi|^2}{\left( \mu^2 + |x-\xi|^2 \right)^N} \\
&=& \mu^{2(2-N)} \omega_N \frac{|\xi|^N}{2^N} \sup_{|x|\, \leq \, |\xi|/2}\frac{C^2_N \,
  |x-\xi|^2}{\left( \mu^2 + |x-\xi|^2 \right)^N} \\
  & \leq& \mu^{2(2-N)} \omega_N \frac{|\xi|^N}{2^N} \sup_{|x|\, \leq \, |\xi|/2}\frac{C^2_N \,
  \left(\,|x|+|\xi| \,\right)^2}{\left( \mu^2 +{\left|\, |x|-|\xi|\,\right|}^2 \right)^N} \\
  & \leq & \frac{9}{4} \omega_N \frac{|\xi|^N}{2^N}\,\frac{C^2_N \,
  |\xi|^2}{\left( \mu^2 +|\xi|^2 /4 \right)^N}
\end{eqnarray*}
where $\omega_N$ is the measure of $S^{N-1}= \left\{x\,\in \, \R^N\,:\, |x|=1 \right\}$. From $N> 4$,  we infer
\begin{equation*}
\frac{ C^2_N \,
  |\xi|^{N+2} }{\left( \mu^2 +|\xi|^2 /4 \right)^N}\, \rightarrow \, 0 \, \quad\text{as $|\xi| \rightarrow +\infty$.}
\end{equation*}
Finally, we deduce
\begin{equation*}
\int_{|x|\, \leq \, |\xi|/2} \left| \nabla_{x}z_0 \left(\frac{x-\xi}{\mu}\right)\right|^2 \, dx\, \rightarrow \, 0 \quad\text{as $\mu\to \bar \mu \in (0,+\infty)$, $|\xi| \to +\infty$}.
\end{equation*}
On the other hand, we have
\begin{equation*}
\int_{|x|\, > \, |\xi|/2} \left| \nabla_{x}z_0 \left(\frac{x-\xi}{\mu}\right)\right|^2 \, dx \, \leq \, \mu^{N-2} \, \int_{|\mu x + \xi |\, > \, |\xi|/2} \left| \nabla_{x} z_0(x) \right|^2 \, dx.
\end{equation*}
Since $| \nabla_{x}z_0 |^2 \,\in \, L^1 (\R^N)$, we deduce that

\begin{equation*}
 \mu^{N-2} \, \int_{|\mu x + \xi |\, > \, |\xi|/2} \left| \nabla_{x} z_0(x) \right|^2 \, dx \, \rightarrow \, 0
\end{equation*}
as $\mu \to \overline{\mu} \in (0, +\infty)$ and $ |\xi| \to + \infty$ and
thus
\begin{equation*}
\int_{|x|\, > \, |\xi|/2} \left| \nabla_{x}z_0 \left(\frac{x-\xi}{\mu}\right)\right|^2 \, dx \, \rightarrow \, 0
\end{equation*}
and
\begin{equation*}
  \alpha_1 \leq
  \mu^{-\frac{(N-2)}{2}}\left\| A \right\|_{L^N (\R^N)}\, \left\|\phi \right\|_{L^{2^*} (\R^N)}\, \left( \int_{\R^N}\left| \nabla_{x}z_0 \left(\frac{x-\xi}{\mu}\right)\right|^2 \, dx \right)^{1/2} \, \rightarrow \, 0
\end{equation*}
as  $\mu \to \overline{\mu} \in (0, +\infty)$ and $ |\xi| \to + \infty$.\\

As regards $\alpha_2$ we know that
\begin{eqnarray*}
\alpha_2 &=& - \re \int_{\R^N}\frac{1}{i}\nabla_{x} \phi (x) \cdot A(x)
e^{-i \sigma} \mu^{-\frac{(N-2)}{2}}z_0 \left(\frac{x-\xi}{\mu}\right) \, dx \\
& \leq & \mu^{-\frac{(N-2)}{2}} \left(\, \int_{\R^N}\left| \, \nabla_{x} \phi (x) \cdot A(x) \right|^{\beta} \, dx\, \right)^{1/ \beta}  \left(\, \int_{\R^N}\left| \,z_0 \left(\frac{x-\xi}{\mu}\right)\, \right|^{2^*} \, dx \right)^{1/2^*} \\
&\leq &  \mu^{-\frac{(N-2)}{2}} \left\|\, \phi\, \right\|_E \, \left\|\,A\, \right\|_{L^N} \left(\, \int_{\R^N}\left| \,z_0 \left(\frac{x-\xi}{\mu}\right)\, \right|^{2^*} \, dx \right)^{1/2^*}
\end{eqnarray*}
with $\beta= 2N/(N+2)$. We notice that
\begin{multline*}
 \int_{\R^N} \left| \,z_0 \left(\frac{x-\xi}{\mu}\right)\, \right|^{2^*} \, dx =
 \int_{| x| \, \leq \,|\xi| /2} \left| \,z_0 \left(\frac{x-\xi}{\mu}\right)\, \right|^{2^*} \, dx \\
{} + \int_{| x| \, > |\xi| /2} \left| \,z_0 \left(\frac{x-\xi}{\mu}\right)\, \right|^{2^*} \, dx.
\end{multline*}
Moreover,
\begin{eqnarray*}
 \int_{| x| \, \leq \,|\xi| /2} \left| \,z_0 \left(\frac{x-\xi}{\mu}\right)\, \right|^{2^*} \, dx & \leq& \omega_N \frac{|\xi|^N}{2^N} \sup_{  |x|\, \leq \, |\xi|/2  } \, \left| z_0 \left(\frac{x-\xi}{\mu}\right)\right|^{2^*} \\
  &=& \omega_N \frac{|\xi|^N}{2^N} \sup_{|x|\, \leq \, |\xi|/2}\mu^{2N}  \kappa^{2^*}_N \frac{ \mu^{2N}  \kappa^{2^*}_N }{\left( \mu^2 + |x-\xi|^2 \right)^N} \\
  & \leq& \mu^{2N} \omega_N \frac{|\xi|^N}{2^N} \sup_{|x|\, \leq \, |\xi|/2}\frac{ \kappa^{2^*}_N}{\left( \mu^2 +{\left|\, |x|-|\xi|\,\right|}^2 \right)^N} \\
  & \leq &   \mu^{2N}   \omega_N \frac{|\xi|^N}{2^N}\,\frac{  \kappa^{2^*}_N       }{\left( \mu^2 +|\xi|^2 /4 \right)^N}
\end{eqnarray*}
where $\omega_N$ is the measure of $S^{N-1}= \left\{x\,\in \, \R^N\,:\, |x|=1 \right\}$. From $N>4$,  we infer
\begin{equation*}
\frac{ \kappa^{2^*}_N \,
  |\xi|^{N} }{\left( \mu^2 +|\xi|^2 /4 \right)^N}\, \rightarrow \, 0 \, \quad\text{as $|\xi| \rightarrow +\infty$.}
\end{equation*}
Finally, we deduce
\begin{equation*}
  \int_{|x|\, \leq \, |\xi|/2} \left|\, z_0 \left(\frac{x-\xi}{\mu}\right)\,\right|^{2^*} \, dx\, \rightarrow \, 0 \,
\end{equation*}
as  $\mu \to \overline{\mu} \in (0, +\infty)$ and $ |\xi| \to + \infty$.

On the other hand, we have
\begin{equation*}
\int_{|x|\, > \, |\xi|/2} \left|\,z_0 \left(\frac{x-\xi}{\mu}\right)\,\right|^{2^*} \, dx \, \leq \, \mu^{N} \, \int_{|\mu x + \xi |\, > \, |\xi|/2} \left|\, z_0(x)\, \right|^{2^*} \, dx.
\end{equation*}
Since $|\,z_0 |^{2^*} \,\in \, L^1 (\R^N)$, we deduce that

\begin{equation*}
 \mu^{N} \, \int_{|\mu x + \xi |\, > \, |\xi|/2} \left|\, z_0(x) \,\right|^{2^*} \, dx \, \rightarrow \, 0
\end{equation*}
as  $\mu \to \overline{\mu} \in (0, +\infty)$ and $ |\xi| \to + \infty$ and thus
\begin{equation*}
\int_{|x|\, > \, |\xi|/2} \left|\,z_0 \left(\frac{x-\xi}{\mu}\right)\,\right|^{2^*} \, dx \, \rightarrow \, 0
\end{equation*}
and
\begin{equation*}
\alpha_2 \, \leq \,
\mu^{-\frac{(N-2)}{2}}\left\| A \right\|_{L^N (\R^N)}\, \left\|\phi \right\|_{E}\, \left( \int_{\R^N}\left| \,z_0 \left(\frac{x-\xi}{\mu}\right)\,\right|^{2^*} \, dx \right)^{1/2^*} \, \rightarrow \, 0
\end{equation*}
as  $\mu \to \overline{\mu} \in (0, +\infty)$ and $ |\xi| \to + \infty$.\\
\\
\\
Conversely, assume that $\mu \, \rightarrow\,+\infty$. Now it is convenient to write
\begin{equation*}
{\left\langle G'_1 (z), \phi \right\rangle}_E = \alpha_1+\alpha_2
\end{equation*}
where
\begin{equation*}
\alpha_1=- \mu \re \int_{\R^N}\frac{e^{i\sigma}}{i} \nabla_y z_0(y) \cdot A(\mu y+\xi) \overline{{\phi^*(y)}}\,dy
\end{equation*}
 and
\begin{equation*}
\alpha_2=- \mu \re \int_{\R^N}\frac{1}{i} \nabla_y {\phi^*(y)}  \cdot A(\mu y+\xi)e^{-i\sigma} z_0(y)   \,dy.
\end{equation*}
The H\"{o}lder inequality implies that
 \begin{equation*}
 \alpha_1 \, \leq \,  \mu  {        \left\|\, \phi^* \,\right\|}_{L^{2^*}} \left(\, \int_{\R^N} \left( \, \nabla_y z_0(y) \cdot A(\mu y+\xi)\right)^\beta \,dy \, \right)^{1/\beta}
 \end{equation*}
 where $1/2^* +1/ \beta=1$ so $\beta= 2N/(N+2)$. By assumptions \textbf{(A1)}, we can fix $r \in (1,(N+2)/2)$ such that $A^\beta \, \in \, L^r (\R^N)$. Moreover, let $s= r /(r-1)$. It is immediate to check that $\beta s\, > \, 2$ and then $|\nabla_y z_0|^{\beta s} \, \in \, L^1(\R^N)$. By \textbf{(A1)} and the H\"{o}lder inequality, we deduce that:
 \begin{multline*}
  \left(\, \int_{\R^N} \left( \, \nabla_y z_0(y) \cdot A(\mu y+\xi)\right)^\beta \,dy \, \right)^{1/\beta}  \\
\leq   \left(\, \int_{\R^N} \left( \, \nabla_y z_0(y)\right)^{\beta s} \,dy \, \right)^{1/\beta s} \left(\, \int_{\R^N} \left( \,A(\mu y+\xi) \right)^{\beta r} \,dy \, \right)^{1/\beta r }\\
\leq \mu^{-\frac{N}{\beta r}} {\left\|\nabla_y z_0(y)\right\|}_{L^{\beta s}}                    \left(\, \int_{\R^N} \left( \,A(\mu y+\xi) \right)^{\beta r} \,dy \, \right)^{1/\beta r }
 \end{multline*}
As a consequence, by the above inequality, we infer for $\mu$ small:
\begin{eqnarray*}
 \alpha_1 &\leq& \mu^{1 -\frac{N}{\beta r}}{\left\|\nabla_y z_0(y)\right\|}_{L^{\beta s}}
  \left(\, \int_{\R^N} \left( \, A(\mu y+\xi) \right)^{\beta r} \,dy \, \right)^{1/\beta r }
  \left\|\, \phi^* \,\right\|_{L^{2^*}} \\
 &\leq& \mu^{1 -\frac{N}{\beta r}} C'_N  {\left\|z_0 \right\|}_E {\left\|  A  \right\|}_{L^{\beta r}}{\left\| \phi^* \right\|}_{E}
\end{eqnarray*}
Analogously,
\begin{eqnarray*}
  \alpha_2 &\leq& \mu {\left(\, \int_{\R^N}{\left( \, \left|\nabla_y \phi^{*}(y) \right| \, \left|A(\mu y+\xi) \right| \, dy \,\right) }^{\beta} \right)}^{1/\beta}
  {\left( \, \int_{\R^N}{\left| z_0(y) \right|}^{2^*} \, dy \,\right)}^{1/2^*}\\
  &\leq& \mu^{1 -\frac{N}{\beta r}} {\left\|\,z_0 \,\right\|}_{L^{2^*}}  \left(\, \int_{\R^N} \left( \, \nabla_y \phi^*(y)\right)^{\beta s} \,dy \, \right)^{1/\beta s}                   \left(\, \int_{\R^N} \left( \,A( y) \right)^{\beta r} \,dy \, \right)^{1/\beta r }\\
  &\leq& \mu^{1 -\frac{N}{\beta r}} C''_N \left( {\left\| \,z_0 \,\right\|}_E \,   {\left\| \,A  \, \right\|}_{L^{\beta r}} \,  {\left\| \,  \phi ^* \,  \right\|}_E    \right)
\end{eqnarray*}
Since $\beta = 2N/(N+2)$, we deduce $1 -\frac{N}{\beta r} < 0$. The
conclusion follows immediately from Lemma \ref{lemma:4.8}.
\end{proof}

\begin{proposition}
  Assume that there exists $\xi \in \R^N$ with $V(\xi) \neq 0$. Then
\begin{equation} \label{eq:47}
\lim_{\mu\,\rightarrow \, 0^+} \frac{\Gamma(\mu, \xi)}{\mu^2} =
\frac{1}{2} V(\xi) \int |z_0|^2.
\end{equation}
In particular, $\Gamma$ is a non-constant map.
\end{proposition}
\begin{proof}
  If $V(\xi)\neq 0$ for some $\xi \,\in\,\R^N$, we can immediately
  check that $\Gamma(\mu, \xi)$ is not identically zero. More
  precisely, we prove that for every $\xi\in\R^N$ there holds
\begin{equation} \label{eq:48}
\lim_{\mu\,\rightarrow \, 0^+} \frac{\Gamma(\mu, \xi)}{\mu^2} =
\frac{1}{2} V(\xi) \int_{\R^N} |z_0|^2.
\end{equation}
Indeed,  after a suitable change of variable,
\begin{eqnarray} \label{eq:14}
\lim_{\mu\,\rightarrow \, 0^+} \frac{G_2(z_{\mu, \xi})}{\mu^2}&=&
\lim_{\mu\,\rightarrow \, 0^+} \frac{1}{2} \int_{\R^N} (|A(\mu y + \xi)|^2
|z_0(y)|^2 + \frac{1}{2} \int_{\R^N} V(\mu y + \xi)  |z_0(y)|^2
 \,dy \notag \\
&=&  \frac{1}{2} |A(\xi)|^2 \int_{\R^N} |z_0(y)|^2\,dy + \frac{1}{2} V(\xi)
\int_{\R^N} |z_0(y)|^2\,dy.
\end{eqnarray}
To complete the proof of \eqref{eq:48},  we need to study $\lim_{\mu\,\rightarrow \, 0^+}\frac{1}{2 \mu^2}
\left\langle G'_1(z_{\mu, \xi}),
\phi_{\mu,\xi}\right\rangle$.

In Lemma \ref{tom}, we have showed that
\[
{\left\langle G'_1( e^{i \sigma }z_{\mu, \xi}), \phi_{\mu,\xi}\right\rangle} = -
{\langle f''_0(z_{\mu, \xi} e^{i \sigma}) \phi_{\mu,\xi},
\phi_{\mu,\xi} \rangle}= - {\langle f''_0(z_{0} e^{i \sigma})
\phi^*_{\mu,\xi}, \phi^*_{\mu,\xi} \rangle}
\]
where $\phi^*_{\mu,\xi}(\frac{x - \xi }{\mu}) = \mu^{N/2 -1}
\phi_{\mu,\xi}(x)$ and $f''_0 (z_0 e^{i \sigma})\phi^{*}_{\mu,\xi}=
k_{\mu,\xi}$, where
\begin{equation*}
k_{\mu,\xi}(y)= \frac{2}{i}\mu e^{i \sigma} \nabla_{y} z_0 (y)
\cdot A(\mu y+\xi) + \frac{\mu^2}{i}e^{i \sigma}
\operatorname{div}_y A(\mu y+\xi) z_0(y).
\end{equation*}
As $\mu \to 0^+$, we have $k_{\mu,\xi} \to k_\xi$, where
\[
k_\xi (x) : = \frac{2}{i} \ e^{i \sigma} \nabla_{y} z_0 (y) \cdot
A(\xi).
\]
Let us define $\psi_\xi(x) = \lim_{\mu \to 0^+}  \frac{L_{z_0}
k_{\mu,\xi}}{\mu}= \lim_{\mu \to 0^+} \frac{
\phi^{*}_{\mu,\xi}}{\mu}$.  We have that
\begin{equation}\label{ixd}
f''_0 (z_0 e^{i \sigma})\psi_{\xi}= \frac{2}{i}\ e^{i \sigma}
\nabla_{x} z_0 (y) \cdot A(\xi).
\end{equation}
Setting $g_\xi (x) = e^{-i \sigma} \psi_\xi (x)$, we have that for any
$\phi \in D^{1,2}(\R^N, \R)$
\[
\langle f''_0 (z_0 e^{i \sigma}) e^{i\sigma} g_\xi ,e^{i\sigma} \phi
\rangle = \re \int \frac{2}{i}\ e^{i \sigma} \nabla_{y}
z_0 (y) \cdot A(\xi)e^{-i\sigma} \phi \, dx =0.
\]
This means that for any $\phi \in {\mathcal D}^{1,2}(\R^N, \R)$
\begin{eqnarray*}
0&=&\langle f''_0 (z_0 e^{i \sigma}) e^{i\sigma} g_\xi ,e^{i\sigma} \phi
\rangle \\
&=& \re \int \nabla (e^{i \sigma} g_\xi) \cdot \overline{\nabla
(e^{i \sigma} \phi)} -
\re \int |z_0|^{2^*-2} e^{i \sigma} g_\xi  \overline{e^{i \sigma} \phi}\\
&-& \re (2^*-2) \int   |z_0|^{2^*-4}
\re (  e^{i \sigma}z_0 e^i {\sigma} g_\xi) e^{i \sigma}z_0 \overline{e^{i \sigma} \phi}\\
&=& \re \int \nabla ( g_\xi) \cdot \overline{\nabla  \phi} - \re \int      |z_0|^{2^*-2}  g_\xi  \overline{  \phi}\\
&-& \re (2^*-2) \int   |z_0|^{2^*-4} \re (z_0  g_\xi) z_0       \overline{  \phi}\\
&=&   \int \nabla (\re  g_\xi) \cdot \overline{\nabla \phi} -  \int      |z_0|^{2^*-2}  \re g_\xi  \overline{  \phi}\\
&-&  (2^*-2) \int   |z_0|^{2^*-4}  \re  g_\xi {z_0}^2       \overline{  \phi}\\
&=& \langle f''_0 (z_0 ) \re g_\xi , \phi
\rangle.
\end{eqnarray*}

It follows that $\re g_\xi =0$ as $\phi_{\mu, \xi} \in
\left( T_{e^{i \sigma}z_{\mu, \xi}}Z \right)^\perp$.
 Therefore $\psi_\xi
(x)= i e^{i\sigma} r_\xi(x)$ with $r_\xi \in D^{1,2}(\R^N, \R)$.
Now we test $(\ref{ixd})$ against  functions of the type $ i e^{i
\sigma} \omega (x)$, $\omega \in D^{1,2}(\R^N, \R)$.

It results:
\begin{eqnarray*}
\re \int  \frac{2}{i}\ e^{i \sigma}
\nabla_{x} z_0 (x) \cdot A(\xi) \overline{i e^{i \sigma} w}&=&\left\langle  f''_0 (z_0 e^{i \sigma})\psi_{\xi}, i e^{i \sigma} w     \right\rangle\\
&=& \re \int \nabla r_\xi \cdot \nabla w -  \re \int      |z_0|^{2^*-2}   r_\xi  w\\
&-&  \re (2^*-2) \int   |z_0|^{2^*-4}  \re (i z_0 r_\xi) z_0  i w
\end{eqnarray*}
or equivalently
\begin{eqnarray*}
\re \int \nabla   r_\xi \cdot \nabla w -  \re \int      |z_0|^{2^*-2}   r_\xi  w =
-\re \int 2
\nabla_{x} z_0 (x) \cdot A(\xi) w.
\end{eqnarray*}
We deduce that $r_\xi$ satisfies the equation
\begin{equation}\label{imps}
- \Delta r_\xi(x) - |z_0|^{2^* -2} r_\xi(x) = - 2 \nabla z_0 \cdot A(\xi).
\end{equation}
We notice that the function $u(x)=z_0(x) A(\xi) \cdot x$ solves the
equation $(\ref{imps})$, as
 $\Delta u = \Delta z_0 A(\xi) \cdot x +
z_0 \Delta (A(\xi) \cdot x) + 2 \nabla z_0 \cdot \nabla(A(\xi) \cdot
x)= \Delta z_0 A(\xi) \cdot x + 2 \nabla z_0 \cdot A(\xi)$.

Since $iz_0(x) (A(\xi) | x) e^{i \sigma}$ belongs to $\left( T_{e^{i
\sigma}z_0}Z \right)^\perp$, we deduce that $\psi_\xi(x) = i e^{i \sigma}
z_0(x) A(\xi) \cdot x$ and thus
\begin{eqnarray*}
 \lim_{\mu\,\rightarrow \, 0^+}\frac{1}{2}
\frac{\left\langle G'_1(z_{\mu, \xi}),
\phi_{\mu,\xi}\right\rangle}{\mu^2} &=& - \re
\int_{\R^N} \frac{1}{i}\ e^{i \sigma} \nabla_{y} z_0 (y) \cdot A(\xi)
\ \overline{i e^{i \sigma}z_0 A(\xi) \cdot x} \, dx \\
&=& \int_{\R^N}  \ \nabla_{y} z_0 (y) \cdot
A(\xi) z_0 \ A(\xi) \cdot x \, dx.
\end{eqnarray*}

Since we have
\[
\int_{\R^N}  \nabla_{y} z_0 (y) \cdot A(\xi) z_0 A(\xi)
\cdot x \, dx
 = -\int_{\R^N} \nabla_{y} z_0 (y) \cdot A(\xi) z_0 A(\xi) \cdot x  \, dx
-\int_{\R^N} |A(\xi)|^2 z_0^2 \, dx,
\]
we conclude that
\begin{equation} \label{eq:15}
 \lim_{\mu\,\rightarrow \, 0^+}\frac{1}{2}
\frac{\left\langle G'_1(z_{\mu, \xi}),
\phi_{\mu,\xi}\right\rangle}{\mu^2} = \int _{\R^N} \ \nabla_{y} z_0 (y)
\cdot A(\xi) z_0 A(\xi) \cdot x \, dx
 = - \frac{1}{2} \int_{\R^N} |A(\xi)|^2 z_0^2 \, dx.
\end{equation}
Therefore we have that
 \[
 \lim_{\mu\,\rightarrow \, 0^+}
\frac{\Gamma(\mu, \xi)}{\mu^2} = \lim_{\mu\,\rightarrow \, 0^+}
\frac{1}{\mu^2}( G_2(\mu, \xi)+ \frac{1}{2}{\left\langle
G'_1(z_{\mu, \xi}), \phi_{\mu,\xi}\right \rangle}) =
 \frac{1}{2} V(\xi) \int_{\R^N} |z_0|^2.
\]
\end{proof}
\begin{remark}
  The presence of a non-trivial potential $V$ is crucial in the
  previous Proposition. Otherwise, from \eqref{eq:14} and
  \eqref{eq:15} we would simply get that $\lim_{\mu \to 0+}
  \frac{\Gamma(\mu,\xi)}{\mu^2} = 0$, and $\Gamma$ might still be a
  constant function. Hence $V$ is in competition with $A$. It would be
  interesting to investigate the case in which $V=0$ identically. We
  conjecture that some additional assumptions on the shape of $A$
  should be made.
\end{remark}

\section{Proof of the main result}

We recall the following abstract theorem from \cite{am99}. See also
\cite{am05}.

\begin{theorem} \label{th:ab} Assume that there exist a set $ A
  \subseteq Z$ with compact closure and $z_0 \in A$ such that
\begin{equation}
\Gamma (z_0) < \inf_{z \in \, \partial A} \Gamma (z)\, \, (\hbox{resp.}\,
\Gamma (z_0) > \sup_{z \in \,  \partial A} \Gamma (z)).\nonumber
\end{equation}
Then, for $\varepsilon$ small enough, $f_ \varepsilon$ has at least a
critical point $u_ \varepsilon \in Z_ \varepsilon$ such that
\begin{equation}
f_0 (z) + \varepsilon^2  \inf_{A} \Gamma + o(\varepsilon^2) \leq f_ \varepsilon (u_ \varepsilon) \leq f_0 (z) +\varepsilon^2 \sup_{\partial A} \Gamma + o (\varepsilon^2) \nonumber
\end{equation}
\begin{equation*}
(\hbox{resp.} \, f_0 (z) + \varepsilon^2 \inf_{\partial A} \Gamma +
o(\varepsilon^2) \leq f_ \varepsilon (u_ \varepsilon) \leq f_0 (z)
+\varepsilon^2 \sup_{A} \Gamma + o(\varepsilon^2)).
\end{equation*}
Furthermore, up to a subsequence, there exists $\overline{z} \in A$
such that $ u_ {\varepsilon_n} \rightarrow \overline{z}$ in $E$ as
$\varepsilon_n \rightarrow 0$.
\end{theorem}

We can finally prove our main existence result for
equation~\eqref{tp}. According to Remark \ref{rem:orbits}, we will use
the term \textit{solution} rather than the more precise $S\sp
1$--\textit{orbit of solutions}.

\begin{theorem}\label{mainresult}
  Retain assumptions \textbf{(N), (A1--2), (V)}. Assume that $V (\xi)
  \neq 0$ for some $\xi \in \R^N$. Then, there exists $\ge_0 >0$ such
  that for all $\ge \in (0,\ge_0)$ equation \eqref{tp} possesses at
  least one solution $u\sb\ge \in E$. If $V$ is a changing sign
  function, then there exists two solutions of equation \eqref{tp}.
\end{theorem}
\begin{proof}
  Under our assumptions, the Melnikov function $\Gamma$, extended
  across the hyperplane $\{\mu =0\}$ by reflection, is not constant
  and possesses at least a critical point (either a minimum or a
  maximum point).  We can therefore invoke Theorem~\ref{th:ab} to
  conclude that there exists at least one solution $u_\ge$ to
  \eqref{tp}, provided $\ge$ is small enough. If there exist points
  $\xi_i \in \R^N$, $i=1,2$, such that $V(\xi_1) V(\xi_2) <0$, then it
  follows from the previous Proposition that $\Gamma$ must change sign
  near $\{ \mu = 0\}$. In particular, it must have \emph{both} a
  minimum \emph{and} a maximum. Hence there exist two different
  solutions to (\ref{tp}).
\end{proof}

\begin{remark}\label{fin}
  Consider equation \eqref{tp}. It is clear that our main theorem
  still applies for any $\alpha \in [1,2)$. Indeed, in the expansion
  \eqref{eq:13bis}, the lowest order term in $\varepsilon$ is
\[
\varepsilon^\alpha \int_{\R^N} V z^2 \, dx,
\]
and consequently the magnetic potential $A$ no longer affects the
finite-dimensional reduction. In some sense, we have treated with the
more all the details the ``worst'' situation in the range $1 \leq
\alpha \leq 2$.
\end{remark}

\section*{Acknowledgement}

The authors would like to thank V.~Felli for some useful discussions
about the proof of Lemma 3.2.



\begin{thebibliography}{99}



\bibitem{ab} A.~Ambrosetti, M.~Badiale, {\it Variational perturbative
    methods and bifurcation of bound states from the essential
    spectrum}, Proc. Royal Soc. Edinburgh {\bf 128 A} (1998),
  1131--1161.

\bibitem{ABC} A.~Ambrosetti, M.~Badiale, S.~Cingolani, {\it
    Semiclassical states of nonlinear {S}chr\"{o}\-dinger equations},
  Arch. Ration.  Mech. Anal. {\bf 140} (1997), 285--300.

\bibitem{aap}A.~Ambrosetti, J.G. Azorero, I. Peral, {\it Perturbation of $\Delta u + u^{(N+2)/(N-2)}=0$,
the Scalar Curvature Problem in $\R^N$, and Related Topics}, Journal
of Functional Analysis {\bf 165} (1999), 117--149.


\bibitem{am99} A.~Ambrosetti, A.~Malchiodi, \textit{A multiplicity
    result for the Yamabe problem on \(S^n\)}, Journal of Functional
  Analysis {\bf 168} (1999), 529--561.

\bibitem{am05} A.~Ambrosetti, A.~Malchiodi, ``Perturbation methods and
  semilinear elliptic problems on \(\mathbb{R}\sp n\)'', Progress in
  Mathematics 240, Birkh\"{a}user Verlag, 2006.


\bibitem{ariolisz} G.~Arioli, A.~Szulkin, {\it A semilinear
    {S}chr\"{o}dinger equations in the presence of a magnetic field},
  Arch. Ration. Mech. Anal.  {\bf 170} (2003), 277--295.

\bibitem{avhesi} J.~Avron, I.~Herbst, B.~Simon, {\it {S}chr\"{o}dinger
    operators with magnetic fields I}, Duke Math. J.  {\bf 45} (1978),
  847--883.

\bibitem{blio} H.\ Berestycki, P.\ L.\ Lions, \textit{Nonlinear scalar field equations
  I and II}, Arch.\ Ration.\ Mech.\ Anal. \textbf{82} (1983) 313--345
  and 347--375.

\bibitem{bs} F. A.~Berezin, M. A.~Shubin, ``The {S}chr\"odinger equation'',
  Mathematics and its Applications (Soviet Series), 66.  Kluwer Academic
  Publishers Group, 1991


\bibitem{bgm} M.~Berger, P.~Gauduchon, E.~Mazet, ``Le spectre d'une vari\'et\'e
  riemannienne'', Lecture Notes in Mathematics, Vol. 194, Springer-Verlag, 1971.

\bibitem{best} A.~Bernoff, P.~Stenberg, {\it Onset of superconductivity in
    decreasing fields or general domains}, J.  Math. Phys. {\bf 39} (1998),
  1272--1284.

\bibitem{caotang} D.~Cao, Z.~Tang, {\it Existence and uniqueness of
    multi--bump bound states of nonlinear Schr\"{o}dinger equations with
    electromagnetic fields}, J. Differential Equations \textbf{222} (2006),
  381--424.

\bibitem{chsz} J.~Chabrowski, A.~Szulkin, {\it On the {S}chr\"odinger
    equation involving a critical {S}obolev exponent and magnetic field},
  Topol. Methods Nonlinear Anal., {\bf 25} (2005), 3--21.

\bibitem{c} S.~Cingolani, {\it Positive solutions to perturbed elliptic
    problems in {$\mathbb{R} \sp N$} involving critical {S}obolev exponent},
  Nonlinear Anal. \textbf{48} (2002), 1165--1178.



\bibitem{ci} S.~Cingolani, {\it Semiclassical stationary states of Nonlinear
    Schr\"{o}dinger equations with an external magnetic field}, J.
  Differential Equations \textbf{188} (2003), 52--79.


\bibitem{cipi} S.~Cingolani, A.~Pistoia,
{\it Nonexistence of single blow-up solutions for a nonlinear
Schrödinger equation involving critical Sobolev exponent},  Z.
Angew. Math. Phys. \textbf{55}   (2004),  201--215.

\bibitem{cs} S.~Cingolani, S.~Secchi, {\it Semiclassical limit for nonlinear
    {S}chr\"{o}dinger equations with electromagnetic fields}, J. Math. Anal.
  Appl.  \textbf{275} (2002), 108--130.

\bibitem{cs1} S.~Cingolani, S.~Secchi, {\it Semiclassical states for
    NLS equations with magnetic potentials having polynomial growths},
  J. Math. Phys. \textbf{46} (2005), 1--19.

\bibitem{cs05} S.~Cingolani, S.~Secchi, {\it Multipeak solutions for NLS
    equations with magnetic fields in semiclassical regime}, to appear.


\bibitem{EL} M.~Esteban, P.L.~Lions, {\it Stationary solutions of
    nonlinear {S}chr\"{o}dinger equations with an external magnetic
    field}, in PDE and Calculus of Variations, in honor of
  E.~De~Giorgi, Birkh\"{a}user, 1990.


\bibitem{hel1} B.~Helffer, {\it On Spectral Theory for Schr\"odinger
Operators with Magnetic Potentials}, Advanced Studies in Pure
Mathematics vol. 23, 113--141 (1994).

\bibitem{RS} B.~Helffer,
{\it Semiclassical analysis for Schr\"odinger operator with magnetic
wells}, in  Quasiclassical methods (J.~Rauch, B.~Simon Eds.). The
IMA Volumes in Mathematics and its applications vol. 95,
Springer--Verlag New--York 1997.


\bibitem{helffermorame} B.~Helffer, A.~Morame, {\it Magnetic bottles
    in connection with superconductivity}, J. Functional Anal. \textbf{93
    A} (2001), 604--680.

\bibitem{ku} K.~Kurata, {\it Existence and semi-classical limit of the
    least energy solution to a nonlinear {S}chr\"{o}dinger equation
    with electromagnetic fields}, Nonlinear Anal. \textbf{41} (2000),
  763--778.


\bibitem{lupan} K.~Lu, X.-B. Pan, {\it Surface nucleation of
    superconductivity in 3-dimensions}, J. Differential Equations \textbf{
    168} (2000), 386--452.


\bibitem{rs} M.~Reed, B.~Simon, ``Methods of Modern Mathematical
  Physics'', vol.II, Academic Press, 1975.


\bibitem{ss} S.~Secchi, M.~Squassina, {\it On the location of spikes
    for the {S}chr\"{o}dinger equations with electromagnetic field},
  Commun. Contemp. Math. \textbf{7} (2005), 251--268.


\bibitem{shen} Z.~Shen, {\it Eigenvalue asymptotics and exponential
    decay of the eigenfunctions for Schr\"odinger operators with
    magnetic fields}, Trans. Amer. Math. Soc.  \textbf{348} (1996),
  4465--4488.

\bibitem{simon} B.~Simon, \textit{Maximal and minimal Schr\"odinger
    forms}, J. Operator Theory {\bf1} (1979), 37--47.

\bibitem{sulem} C.~Sulem, P.L.~Sulem, ``The Nonlinear Schr\"odinger
  Equation'', Self-Focusing and Wave Collapse, Springer 1999.

\bibitem{ta} G.~Talenti, \textit{Best constant in {S}obolev
    inequality}, Ann. Mat. Pura Appl., {\bf 110} (1976), 353--372.

\end{thebibliography}
\end{document}